\documentclass{article}

\usepackage{arxiv}
\usepackage[utf8]{inputenc} 
\usepackage[T1]{fontenc}    
\usepackage{url}            
\usepackage{booktabs}       
\usepackage{amsfonts}       
\usepackage{microtype}      
\usepackage{lipsum}
\usepackage{amssymb}
\usepackage{amsmath}
\usepackage{amscd}
\usepackage{graphicx}
\usepackage{graphics}
\usepackage[noadjust]{cite}
\usepackage{amsthm}
\usepackage{caption}
\usepackage{multirow}
\usepackage{array}
\usepackage{rotating}
\usepackage[norelsize,ruled,linesnumbered]{algorithm2e}
\usepackage{algorithm2e,setspace}
\usepackage{hyperref}

\usepackage{indentfirst}
\usepackage{tikz}
\usepackage{calc}
\usepackage{epsfig}
\usepackage{natbib}
\usepackage[makeroom]{cancel}
\usepackage{subcaption}

\title{Implementing bound constraints and total-variation regularization in extended full waveform inversion with the alternating direction 
method of multiplier: application to large contrast media}

\author{
  Hossein S.~Aghamiry\\
  Institute of Geophysics, University of Tehran, Tehran, Iran.
  \texttt{h.aghamiry@ut.ac.ir} \\
  Geoazur, Universit\'e C\^ote d'Azur, CNRS, IRD, OCA, Valbonne, France. 
  \texttt{aghamiry@geoazur.unice.fr} 
   \And
 Ali Gholami \\
  Institute of Geophysics, University of Tehran, Tehran, Iran.
  \texttt{agholami@ut.ac.ir} 
  \And
  St\'ephane Operto \\
  Geoazur, Universit\'e C\^ote d'Azur, CNRS, IRD, OCA, Valbonne, France. 
  \texttt{operto@geoazur.unice.fr}
  }

\begin{document}
\maketitle

\begin{abstract}
Full waveform inversion (FWI) is a waveform matching procedure, which can provide a subsurface model with a wavelength-scale resolution. However, this high resolution makes FWI prone to cycle skipping, which drives the inversion to a local minimum when the initial model is not accurate enough. Other sources of nonlinearities and ill-posedness are noise, uneven illumination, approximate wave physics and parameter cross-talks. All these sources of error require robust and versatile regularized optimization approaches to mitigate their imprint on FWI while preserving its intrinsic resolution power. To achieve this goal, we implement bound constraints and total variation (TV) regularization in the so-called frequency-domain wavefield-reconstruction inversion (WRI) with the alternating direction method of multipliers (ADMM). In the ADMM framework, WRI relies on an augmented Lagrangian function, a combination of penalty and Lagrangian functions, to extend the FWI search space by relaxing the wave-equation constraint during early iterations. Moreover, ADMM breaks down the joint wavefield reconstruction plus parameter estimation problem into a sequence of two linear subproblems, whose solutions are coordinated to provide the solution of the global problem. The decomposability of ADMM is further exploited to interface in a straightforward way bound constraints and TV regularization with WRI via variable splitting and proximal operators. The resilience of our regularized WRI formulation to cycle skipping and noise as well as its resolution power are illustrated with two targets of the large-contrast BP salt model. Starting from a 3Hz frequency and a crude initial model, the extended search space allows for the reconstruction of the salt and subsalt structures with a high fidelity. The TV regularization filters out the imprint of ambient noise and artifacts associated with multi-scattering and Gibbs effects, while fostering large-contrast reconstruction. Compared to other TV-regularized WRI implementations, the proposed method is easy to tune due to its moderate sensitivity to penalty parameters and does not require a prior guess of the TV-norm ball.
\end{abstract}

\keywords{ FWI, WRI, TV regularization, Bound constraints, ADMM, Split Bregman.}

\section{Introduction}
%
%
%
During the last decade, full waveform inversion (FWI) has been used to estimate subsurface parameters (P and S wavespeeds, density, attenuation, anisotropic parameters) with a resolution close to the seismic wavelength by matching recorded and synthetic seismograms \citep{Tarantola_1984_ISR,Pratt_1998_GNF,Virieux_2009_OFW}. From the numerical optimization viewpoint, the data-fitting/parameter-estimation problem underlying FWI is a nonlinear partial differential equation (PDE)-constrained optimization problem, where the equality constraint is the wave equation and the optimization parameters are embedded in the coefficients of the PDE. Due to the computational burden of multiple source modelling and the size of the data and parameter spaces, this PDE-constrained optimization problem is solved with iterative local (linearized) optimization techniques, namely gradient-based methods \citep{Nocedal_2006_NOO}. Moreover, it is often solved with a reduced-space formulation, which means that the full search space that encompasses the unknown wavefield and the subsurface parameters is first projected onto the parameter space by computing exactly the incident wavefields in the current subsurface model before updating this later \citep{Haber_2000_OTS,Askan_2007_FWI,Epanomeritakis_2008_NCG}. 
It is well acknowledged that the oscillating nature of seismic signals makes the reduced-space formulation highly nonlinear as the modelled seismograms computed in the current subsurface model may be too far away from the recorded ones to satisfy the cycle-skipping criterion, that is the modelled seismograms should predict the recorded traveltimes with an error lower than half a period \citep[e.g.][]{Virieux_2009_OFW}. 

%
%

Beyond cycle skipping, other sources of nonlinearity and ill-posedness such as noise, uneven subsurface illumination, approximate wave physics and parameter cross-talks in multi-parameter reconstruction require the use of stabilizing or regularization techniques that drive the inversion towards subsurface models that satisfy some a priori assumptions. 
Among the penalization techniques, Tikhonov regularization is probably the most popular one and seeks to penalize the roughness of the subsurface model to force smooth reconstruction (see \citet{Benning_2018_MRM} for a review). As the subsurface may be better represented by piecewise smooth media with potentially sharp contrasts as in presence of salt, edge-preserving techniques such as total-variation (TV) regularization have been proposed to steer the inversion to the space of blocky structured models. TV regularization has been applied on several geophysical applications such as FWI \citep{Askan_2007_FWI,Anagaw_2011_RAF,Guitton_2012_BRS,Maharramov_2015_RST,Peters_2016_CVP,Brandsberg-Dahl_2017_FMU,Esser_2018_TVR,Kalita_2018_IBC}, seismic tomography \citep{Gholami_2010_RLN,Loris_2012_IAT}, impedance inversion \citep{Gholami_2015_NMI,Gholami_2016_FAM}, amplitude versus offset (AVO) inversion \citep{Gholami_2017_CNA}, and seismic deconvolution \citep{Gholami_2013_FBS}. For FWI applications, TV regularization can be implemented as a penalty function \citep{Askan_2007_FWI,Anagaw_2011_RAF,Brandsberg-Dahl_2017_FMU,Kazei_2018_SBI} or as a constraint \citep{Peters_2016_CVP,Esser_2018_TVR}. Choosing the most suitable implementation strategy may depend on the prior information on the TV norm of the model and on the optimization method that is used to minimize the objective function \citep{Alkhalifah_2018_FMW}. If the information about the value of TV norm of the model is available, TV regularization can be implemented as a constraint in the FWI objective function. Otherwise, one may resort to a penalty method with the difficulty to design an adaptive penalty parameter, which optimally balances over iterations the relative weight of the data misfit and the total variation of the model in the objective.\\

%
%
%
In this context, the objective of this study is to present a novel implementation of TV regularization and bound-constraints in frequency-domain FWI based upon wavefield reconstruction inversion (WRI).
WRI has been originally proposed by \citet{VanLeeuwen_2013_MLM} to extend the search space and mitigate the risk of cycle skipping accordingly. WRI recasts the PDE-constrained optimization problem underlying FWI into an unconstrained quadratic penalty method, where the penalty term is the 
$\ell_2$ norm of the source residuals (namely, the PDE-constraint violation) that is weighted by a positive penalty parameter $\lambda$. The penalty method relaxes the wave-equation constraint at the benefit of the data fitting during early iterations, hence mitigating the risk of cycle skipping.
To make WRI computationally tractable, \citet{VanLeeuwen_2013_MLM} perform the wavefield reconstruction and the subsurface parameter estimation in an alternating way: first, keeping the current subsurface model fixed, the wavefields, which best jointly fit the observations and satisfy the wave equation in a least-squares sense, are reconstructed for each source; Second, keeping the previously-reconstructed wavefields fixed, the subsurface parameters are estimated by  least-squares  minimization of the source residuals the wave-equation relaxation generated. This cycle being iterated until convergence. A nice property of the alternating-direction strategy is to linearise the parameter-estimation subproblem around the reconstructed wavefield because the wave equation constraint is bilinear.  
However, a significant pitfall of WRI resides in the tuning of the penalty parameter $\lambda$. Ideally, increasing values should be used during iterations to  progressively enforce the wave-equation constraint and, hence satisfy the first-order optimality conditions of the original constrained problem with acceptable precision at the minimizer. A significant issue is that this continuation approach is tedious to implement and the Hessian is ill conditioned for large $\lambda$. Therefore, \citet{VanLeeuwen_2013_MLM} implement WRI with a small preset value of $\lambda$, which leads to slow convergence and a subsurface model of limited accuracy. \\
Later, \citet{vanLeeuwen_2016_PMP} reformulated WRI as a reduced penalty method implemented with a variable projection approach: the closed-form expression of the extended-domain reconstructed wavefield is injected as a function of the subsurface parameters in the penalty function instead of using this wavefield as a passive variable (i.e., independent to the subsurface parameters). Unlike the alternating-direction approach, this variable projection leaves the parameter-estimation subproblem non linear. \citet{vanLeeuwen_2016_PMP} assess their method with a Gauss-Newton method (by opposition to the full Newton counterpart) to mitigate the computational burden. Moreover, using a sparse approximation of the Gauss-Newton Hessian makes the descent direction of the reduced approach identical to that of the alternating-direction WRI of \citet{VanLeeuwen_2013_MLM}. Also, \citet{Aravkin_2017_EQP} analysed the convergence properties of the reduced penalty method when the full Hessian is taken into account and concluded that the variable projection penalty method is insensitive to the penalty parameter. However, this convergence property still needs to be verified against realistic numerical experiments. \\

To make the alternating-direction WRI of \citet{VanLeeuwen_2013_MLM} (referred to as WRI in the following for sake of brevity) more independent to the penalty parameter, \citet{Aghamiry_2019_IWR} have replaced the penalty method by an augmented Lagrangian method \citep{Nocedal_2006_NOO}, leading to the so-called iteratively-refined (IR)-WRI method.  
As in WRI, IR-WRI performs the primal wavefield and parameter updates in an alternating mode, while the Lagrange multipliers (i.e., the dual variables) are updated with a gradient ascent method. As above mentioned, this alternating direction strategy makes the parameter estimation subproblem linear due to the bilinearity of the wave equation constraint. It follows from this linearization that the alternating direction strategy combined with the augmented Lagrangian method is equivalent to an extension of the alternating direction method of multiplier (ADMM) to biconvex problem \citep{Boyd_2011_DOS}. Also, using a scaled form of the augmented Lagrangian, we recast IR-WRI as a penalty method where the right-hand sides (the data and the sources) in the objectives of the penalty function are iteratively updated with the running sum of the data and source residuals (the dual gradient steps). This reformulation of IR-WRI as a penalty method with right-hand side updating clearly draws some similarities and differences with WRI.
The right-hand side updating makes IR-WRI largely insensitive to the penalty parameter for a wide range of preset values \citep[][ Their Figures 2 and 3]{Aghamiry_2019_IWR}. Using a moderate value of the penalty parameter allows for significant wave equation error and improved data fitting during early iterations for search space extension, without preventing the fulfilment of the wave equation constraint with small error at the minimizer \citep[][Chapter 17, Theorem 17.6]{Nocedal_2006_NOO}. This adaptivity makes IR-WRI resilient to cycle skipping as WRI with however a much faster convergence toward a more accurate minimizer. The reader is referred to \citet{Aghamiry_2019_IWR} for a thorough comparative convergence and accuracy analysis of WRI and IR-WRI based upon toy and complex large-contrast synthetic examples.

The objective of this study is to show how to interface TV regularization and bound constraints (hereafter, we refer to it as BTV regularization) with IR-WRI by taking advantage of the alternating-direction strategy implemented in ADMM  and the split-Bregman variable splitting scheme developed by \citet{Goldstein_2009_SBM}. 
More precisely, we recast the BTV regularized IR-WRI as a TV minimization problem subject to constraints, that are the modelled wavefield fit the observables and satisfy the wave equation with prescribed errors, and the model parameters preset bounds. As in WRI and IR-WRI, we solve the wavefield and subsurface parameter subproblems in an alternating mode. However, the later one involves now a combination of $\ell_1$ and $\ell_2$ norms related to the TV minimization and wave-equation error minimization, respectively, with additional bound constraints. 
This is managed by the split-Bregman variable splitting scheme, which de-couples the  $\ell_1$ and $\ell_2$ components and bound constraints of the functional through the introduction of auxiliary variables and solves each related subproblem in sequence \citep{Goldstein_2009_SBM}.\\

We first apply our method on a toy example corresponding to a high-velocity box-shape anomaly embedded in a background model where the velocity increases with depth. Then, we consider two more realistic  examples corresponding to the left and central parts of the large-contrast 2004 BP salt model \citep{Billette_2004_BPB}. We show that the BTV regularized IR-WRI converges to accurate minimizers when we start from a crude initial model and a realistic 3~Hz frequency. We also compare the results of WRI and IR-WRI without any priors, and when performed only with bound constraints and with BTV regularization to highlight the impact of each ingredient upon the quality of the results and the computational burden.  We also assess the resilience of the method to noise by comparing the results that are obtained with noiseless and noisy data.

%
%

This paper is organized as follows. In the first section, we review the principles of WRI and IR-WRI. We first recast FWI as a feasibility problem and review different approaches that are suitable to solve PDE-constrained optimization problems such as penalty and augmented Lagrangian methods. Then, we review how we can easily interface some stabilizing terms with the feasibility problem through variable splitting and ADMM. In the second part, we present the results of the synthetic examples involving the inclusion model with two different starting models and the two targets of the BP2004 model with noiseless and noisy data. The results confirm that the combined use of TV regularization and bound constraint in the ADMM-based IR-WRI method defines a suitable framework to make high-resolution FWI immune to cycle skipping in large-contrast media. 
%
%
%
%
%
\section{Method}
In the following, we first recast frequency-domain FWI as a bi-convex feasibility problem, which can be formulated as a  constrained optimization problem with identically-zero objective function \citep[their appendix A]{Aghamiry_2019_IWR}. Then, we review penalty and augmented Lagrangian methods as optimization techniques to solve this constrained optimization problem with extended search-space, leading to WRI and IR-WRI, respectively. Finally, we interface bound constraints and isotropic TV regularization with IR-WRI by replacing the identically-zero objective function by the TV norm of the subsurface model, and we show how to solve efficiently the regularized IR-WRI with ADMM and split Bregman iterations.
\subsection{WRI and IR-WRI principles}
FWI can be formulated in the frequency domain as the following bi-convex feasibility problem \citep{Aghamiry_2019_IWR}:
\begin{subequations}
\label{feasibility}
\begin{eqnarray}
&& {\text{Find}} ~~~~~~~~~\quad \bold{m} ~\text{and}~ \bold{u} \label{feasibilitya} \\
&& \text{subject to} ~~~~\bold{F(m)u}=\bold{s}\label{feasibilityb}
\end{eqnarray} 
\end{subequations} 
with
\begin{equation} \label{OP}
 \bold{s} =
 \begin{bmatrix}
\bold{d} \\
\bold{b}
\end{bmatrix} 
,~~~~~~~~
 \bold{F}(\bold{m}) =
 \begin{bmatrix}
\bold{P} \\
\bold{A(m)}
\end{bmatrix}, 
\end{equation}
where $\bold{m} \in \mathbb{R}^{N\times 1}$ denotes the vector of discrete model parameters (here, the squared slowness),
$\bold{u} \in \mathbb{C}^{N\times 1}$ the wavefield,
$\bold{b} \in \mathbb{C}^{N\times 1}$ the source term,
$\bold{d} \in \mathbb{C}^{M\times 1}$ the recorded wavefield (data) at receiver locations, and $\bold{P} \in \mathbb{R}^{M\times N}$ is a linear observation operator that samples the modelled wavefield at the receiver positions.  
The matrix $\bold{A(m)} \in \mathbb{C}^{N\times N}$ represents the discretized PDE Helmholtz operator \citep{Pratt_1998_GNF,Plessix_2007_HIS, Chen_2013_OFD}. 
\begin{eqnarray}  
\bold{A(m)} &= & \bold{\Delta} + \omega^2 \bold{C}(\bold{m}) \text{diag}(\bold{m})\bold{B},
\label{helmholtz}
\end{eqnarray}
where $\omega$ is the angular frequency and $\bold{\Delta}$ is the discretized Laplace operator. 
The operator $\bold{C}$ encloses boundary conditions, which can be a function of $\bold{m}$ \citep[e.g., Robin paraxial conditions,][] {Engquist_1977_ABC} or independent from $\bold{m}$ \citep[e.g., sponge-like absorbing boundary conditions such as perfectly-matched layers,][] {Berenger_1994_PML}. Also, the linear operator $\bold{B}$ can be used to spread the "mass" term $\omega^2\bold{C}(\bold{m}) \text{diag}(\bold{m})$ 
over all the  coefficients of the stencil to improve its accuracy following an anti-lumped mass strategy \citep{Marfurt_1984_AFF,Jo_1996_OPF,Hustedt_2004_MGS}.   
In the feasibility problem \ref{feasibility}, we just want to find $\bold{m}$ and $\bold{u}$ that satisfy the constraint (we assume that the feasible set is non-empty, namely the constraint is consistent).
The feasibility problem can be formulated as the following constrained optimization problem with identically-zero objective function \citep[][ Page 128]{Boyd_2004_COO}
\begin{subequations}
\label{feas_opt}
\begin{eqnarray}
&& \min_{\bold{u,m}} ~~~~~~~~~\quad \bold{0} \label{feas_opta} \\
&& \text{subject to} ~~~~\bold{F(m)u}=\bold{s}\label{feas_optb}.
\end{eqnarray} 
\end{subequations} 
\subsubsection{WRI penalty method}
WRI implements this constrained optimization problem with a penalty method.
\begin{equation} \label{leeuwen_FWI_closed}
\min_{\bold{u,m}} ~~ \bold{0} + \|\bold{F(m)u}-\bold{s}\|_{\bold{\Lambda}}^2,
\end{equation}
where $\|\bold{x}\|_{\bold{Q}}^2:= \bold{x}^T\bold{Qx}$ for a vector $\bold{x}$ and a square matrix $\bold{Q}$ with superscript $T$ denoting matrix transposition.
In equation \ref{leeuwen_FWI_closed}, $\bold{\Lambda}$ is a diagonal matrix which includes the penalty parameters $\lambda_0$, $\lambda_1 >0$  on its main diagonal,

\begin{equation} \label{Lambda}
\bold{\Lambda} =
 \left[
    \begin{array}{r@{}c|c@{}l}
  &    \begin{smallmatrix}
       \lambda_0 & & 0 \\
          &\ddots&\\
        0 & & \lambda_0 \rule[-1ex]{0pt}{2ex}
      \end{smallmatrix} 
  & \mbox{\huge0} 
  & \rlap{\kern5mm$M$}\\
      \hline
  &    \mbox{\huge0} &  
       \begin{smallmatrix}\rule{0pt}{2ex}
        \lambda_1 & & 0 \\
          &\ddots&\\
        0 & & \lambda_1
      \end{smallmatrix}    &  \rlap{\kern5mm $N$}
    \end{array} 
\right].
\end{equation}
The objective function in equation \ref{leeuwen_FWI_closed} is a compact writing of the penalty function of \citet{VanLeeuwen_2013_MLM} 
\begin{equation} \label{leeuwen_FWI}
\min_{\bold{u,m}}  \|\bold{Pu}-\bold{d}\|_2^2 + \lambda \|\bold{A(m)u}-\bold{b}\|_2^2,
\end{equation}
where $\lambda = \lambda_1/\lambda_0$. They solve this biconvex minimization problem with an alternating-direction approach to break down the full problem into a sequence of two linear sub-problems:  
A cycle of the algorithm first reconstructs, for each source, the wavefield  $\bold{u}$ that best fits the data and satisfies the wave equation in a least-square senses for the current subsurface model. Then, the subsurface model is updated by minimization of the source residuals with a Gauss-Newton algorithm keeping the reconstructed wavefields fixed. 
A difficulty with the penalty method given by equation \ref{leeuwen_FWI} resides in the tuning of the penalty parameter $\lambda$ during iterations. An increasing values of $\lambda$ should be used during iterations, known as the penalty algorithm, to progressively enforce the wave-equation constraint in iterations and hence satisfy the Karush-Kuhn-Tucker (KKT) optimality conditions \citep{Nocedal_2006_NOO} associated with the original constrained problem with acceptable precision. The main problem is that this continuation strategy is tedious to implement and a large $\lambda$ makes the problem severely ill-conditioned. \\
\subsubsection{IR-WRI augmented Lagrangian method}
To bypass this difficulty, IR-WRI implements the original constrained problem, equation \ref{feas_opt}, 
with the augmented Lagrangian (AL) method \citep{Hestenes_1969_MUL,Nocedal_2006_NOO,Boyd_2011_DOS,Bertsekas_2016_HMC}. 
\begin{equation} \label{AL}
\min_{\bold{u,m}} \max_\bold{v} ~~ \bold{0} + \bold{v}^T[\bold{F(m)u}-\bold{s}]  + \frac{1}{2}\|\bold{F(m)u}-\bold{s}\|_{\bold{\Lambda}}^2,
\end{equation}
where $\bold{\Lambda}$ is defined as in equation \ref{Lambda} and $\bold{v} \in \mathbb{C}^{(M+N)\times 1}$ is the Lagrangian multiplier (known as dual variable).
Comparing the penalty function, equation~\ref{leeuwen_FWI_closed}, and the augmented Lagrangian function, equation~\ref{AL}, clearly shows that the augmented Lagrangian method combines a penalty method with a Lagrangian method. A first advantage of the augmented Lagrangian method relative to the penalty method is to prevent ill-conditioning by introducing explicit estimate of the Lagrange multiplier in the optimization \citep[][ chapter 17]{Nocedal_2006_NOO}. Moreover, the Lagrange multiplier gives the augmented Lagrangian method one more way of improving the accuracy of the minimizer in addition to the penalty parameter, hence allowing for a fixed value to be used for this latter \citep[][ Theorem 17.6]{Nocedal_2006_NOO}.
Applying the alternating direction strategy of WRI on the augmented Lagrangian method leads to an adaptation of the alternating-direction method of multiplier (ADMM) to biconvex problem \citep[][ Section 9.2]{Boyd_2011_DOS}. 
One ADMM iteration first minimizes the augmented Lagrangian function with respect to the primal variables $\bold{u}$ and $\bold{m}$ via a single Gauss-Seidel like iteration (namely, fix one variable and solve for the other) and then update the Lagrangian multiplier via a gradient ascent method.

In the following section, we review each step of the ADMM-based IR-WRI algorithm when equipped with TV regularization and bound constraints.  The reader is also referred to \citet{Aghamiry_2019_IWR} for the detailed IR-WRI algorithm when no regularization is used.

\subsection{BTV-regularized ADMM-based IR-WRI}
To implement TV regularization and bound constraints in IR-WRI, we first recast FWI as a constrained TV minimization problem given by
\begin{subequations}
\label{main1_12tv}
\begin{eqnarray}
&& \min_{\bold{u},\bold{m} \in \mathcal{C}} ~~~~~~~~~\quad \|\bold{m}\|_\text{TV} \\
&& \text{subject to} ~~~~\bold{F(m)u}=\bold{s}, \label{constv}
\end{eqnarray} 
\end{subequations}
where $\bold{F}$ and $\bold{s}$ are defined as in equation \ref{OP}, $\|\bold{m}\|_\text{TV} = \sum  \sqrt{|{\bold{\nabla}}_1 \bold{m}|^2 + |{\bold{\nabla}}_2\bold{m}|^2} $ is the blockiness-promoting isotropic TV norm \citep{Rudin_1992_NTV}, and $\bold{\nabla}_1$ and $\bold{\nabla}_2$ are first-order finite-difference operators in the horizontal and vertical directions, respectively. With notation abuse, the absolute sign, square power, and the square root operations are done component-wise, and the sum runs over all elements (the domain of parameters). 
Also $\mathcal{C} = \{\bold{x} \in \mathbb{R}^{N\times 1}~\vert~ \bold{m}_{{l}} \leq \bold{x} \leq \bold{m}_{{u}}\}$ is the set of all feasible models bounded by the lower bound $\bold{m}_{{l}}$ and the upper bound $\bold{m}_{{u}}$. \\
Compared to the FWI definition given in the previous section, equation~\ref{feas_opt}, we have replaced the identically-zero objective by the TV norm of the model and restricted the space of feasible models to $\mathcal{C}$.
%
%
%
%
Accordingly, the augmented Lagrangian function for the problem defined by equation \ref{main1_12tv} is 
\begin{equation} \label{AL1}
\mathcal{L}_A(\bold{m},\bold{u},\bold{v}) = \|\bold{m}\|_\text{TV} +  \bold{v}^T[\bold{F(m)u}-\bold{s}]  + \frac{1}{2}\|\bold{F(m)u}-\bold{s}\|_{\bold{\Lambda}}^2,
\end{equation}
where the same notations as those of equation~\ref{AL} are used. \\
Equation \ref{AL1} can also be written in a more compact form as
\begin{equation} \label{new_scale}
\mathcal{L}_A(\bold{m},\bold{u},\bold{\bar{s}}) = \|\bold{m}\|_\text{TV} + \frac{1}{2}  \|\bold{F(m)u} - \bold{s}-\bold{\bar{s}}\|_{\bold{\Lambda}}^2-\frac{1}{2} \|\bold{\bar{s}}\|_{\Lambda}^2,
\end{equation}
where $\bold{\bar{s}}=-\bold{\Lambda}^{-1}\bold{v}$ is the scaled dual variable and equation \ref{new_scale} is the scaled form of the augmented Lagrangian (\citep[][ Page 15]{Boyd_2011_DOS} and Appendix A). 
The scaled form recasts the augmented Lagrangian method as a quadratic penalty method where the right-hand sides are updated with the scaled dual variables. 
This highlights similarities and differences between WRI and IR-WRI, since this right-hand side updating is lacking in the former. \\
The method of multipliers seeks to find the saddle point of the scaled augmented Lagrangian \ref{new_scale} through a primal descent - dual ascent updating resulting in the following iteration:
\begin{subequations}
\label{scale_main1}
\begin{eqnarray}
\bold{m}^{k+1},\bold{u}^{k+1} &=& \underset{\bold{u},\bold{m}\in \mathcal{C}}{\arg\min} ~~ 
 \|\bold{m}\|_\text{TV} + \frac{1}{2}  \|\bold{F(m)u} - \bold{s}-\bold{\bar{s}}^k\|_{\bold{\Lambda}}^2 \label{scale_main1a} \\
\bold{\bar{s}}^{k+1} &=& \bold{\bar{s}}^{k} +  \bold{s}- \bold{F}(\bold{m}^{k+1})\bold{u}^{k+1} , \label{scale_main1b}
\end{eqnarray} 
\end{subequations}
for $k=0,1,...$ beginning with a prior estimate $\bold{\bar{s}}^0=\bold{0}$.  The iteration \ref{scale_main1} can be viewed as follows:
we begin with a prior estimate of the dual $\bold{\bar{s}}^0$, and minimize the objective function with respect to the primal variables $\bold{m}$ and $\bold{u}$, equation \ref{scale_main1a}. Subsequently, we maximize the objective function  with respect to the dual variable $\bold{\bar{s}}$ with a gradient ascent method when $\bold{m}$ and $\bold{u}$ are kept fixed, equation \ref{scale_main1b}. The steepest-ascent step, equation \ref{scale_main1b}, shows that the scaled dual variable $\bold{\bar{s}}$ is updated with the residual constraint violation of the current iteration. Remembering that the constraint gathers the observation equation $\bold{Pu}= \bold{d}$ and the wave equation $A(\bold{m}) \bold{u} = \bold{b}$, equation~\ref{OP}, the scaled dual variable $\bold{\bar{s}}$ updates the right-hand sides of the quadratic penalty function, equation~\ref{scale_main1a}, with the running sum of the data and source residuals in iterations. This right-hand side updating describes the well-known iterative solution refinement procedure for ill-posed linear inverse problems as reviewed by \citet[][ Their appendix B]{Aghamiry_2019_IWR}.
This process is iterated until convergence, i.e., when $\bold{F}(\bold{m}^{k+1})\bold{u}^{k+1}=\bold{s}$. In the following, we remove the bar of $\bold{\bar{s}}^k$ for the sake of simplicity.\\
Solving the subproblem \ref{scale_main1a} jointly for the primal variables ($\bold{m}$, $\bold{u}$) is computationally too intensive. 
A splitting method is useful here to break down the joint optimization over $\bold{m}$ and $\bold{u}$ into two subproblems (the readers can refer to \citet{Glowinski_2017_SMC} for an overview of splitting methods). The alternating-direction method of multipliers (ADMM) \citep{Boyd_2011_DOS} 
provides a simple framework to achieve this goal via a single Gauss-Seidel like iteration leading to the following iteration:
\begin{subequations}
\label{main30}
\begin{eqnarray}
\bold{u}^{k+1} &= & \underset{\bold{u}}{\arg\min} ~~ \frac{1}{2}  \|\bold{F}(\bold{m}^k)\bold{u}-\bold{s}-\bold{s}^k\|_{\bold{\Lambda}}^2 \label{main3a0} \\
\bold{m}^{k+1} &= & \underset{\bold{m}\in \mathcal{C}}{\arg\min} ~~ 
 \|\bold{m}\|_\text{TV}+\frac{1}{2}  \|\bold{F(m)u}^{k+1}-\bold{s}-\bold{s}^k\|_{\bold{\Lambda}}^2 \label{main3b0} \\
\bold{s}^{k+1} &=& \bold{s}^{k} +  \bold{s}- \bold{F}(\bold{m}^{k+1})\bold{u}^{k+1} .       \label{main3d00}
\end{eqnarray} 
\end{subequations} 
The ADMM iteration has decomposed the full problem into two subproblems associated with primal variables $\bold{u}$ and $\bold{m}$ by passing the primal update of one subproblem as a passive variable for the next subproblem. We show below that, taking advantage of the bi-convexity of the problem, this alternating-direction approach has linearized the primal subproblem for $\bold{m}$ around the reconstructed wavefield $\bold{u}$. 
The fact that the primal update of one subproblem is passed to the next subproblem implies obviously that the two subproblems are solved in sequence rather than in parallel as in ADMM for linear separable problems. \\
A second modification related to ADMM resides in the updating of the dual variable.
In ADMM, the dual variables are updated only once per iteration after the primal-variable updates, equation~\ref{main3d00}.    
A variant of ADMM, referred to as the Peaceman-Rachford splitting method (PRSM) \citep{peaceman_1955_PRA}, consists of updating the Lagrange multipliers several times, once after the update of each primal variable (see \citet[][ Compare equations 1.3 and 1.4]{He_2014_SPR}.
One issue with PRSM relative to ADMM is that PRSM requires more restrictive assumptions to ensure its convergence, while it is always faster than ADMM whenever it is convergent \citep{He_2014_SPR} .
This issue prompted \citet{He_2014_SPR} to implement a relaxation factor (or, step length) $\alpha\in (0,1)$ to guarantee the strict contraction of the PRSM iterative
sequence. 
Applying the strictly contractive PRSM algorithm to our minimization problem gives
\begin{subequations}
\label{main4}
\begin{eqnarray}
\bold{u}^{k+1} &= & \underset{\bold{u}}{\arg\min} ~~ \frac{1}{2}  \|\bold{F}(\bold{m}^k)\bold{u}-\bold{s}-\bold{s}^k\|_{\bold{\Lambda}}^2 \label{main4a0} \\
\bold{s}^{k+\frac{1}{2}} &=& \bold{s}^{k} +  \alpha[\bold{s} - \bold{F}(\bold{m}^{k})\bold{u}^{k+1}] \label{main4b00}\\
\bold{m}^{k+1} &= & \underset{\bold{m}\in \mathcal{C}}{\arg\min} ~~ 
 \|\bold{m}\|_\text{TV}+\frac{1}{2}  \|\bold{F(m)u}^{k+1}-\bold{s}-\bold{s}^{k+\frac{1}{2}}\|_{\bold{\Lambda}}^2 \label{main4b0} \\
\bold{s}^{k+1} &=& \bold{s}^{k+\frac{1}{2}} +  \alpha[\bold{s}-\bold{F}(\bold{m}^{k+1})\bold{u}^{k+1}] \label{main4b1}.
\end{eqnarray} 
\end{subequations}
The reader is referred to \citet[][ Their Figure 4]{Aghamiry_2019_IWR} for a comparative numerical analysis of the convergence speed of ADMM and PRSM in IR-WRI. In
our numerical tests we found that $\alpha=0.5$ can serve as a suitable value.\\
We now provide the closed-form solution of the two primal subproblems associated with $\bold{u}$ and $\bold{m}$.
\subsubsection{Solving for $\bold{u}$}
The primal subproblem associated with $\bold{u}$, equation \ref{main4a0}, is a linear optimization problem whose solution satisfies in a least-squares sense the following system of linear equations:
\begin{equation}  \label{solve_sub_uu}
\begin{bmatrix}
\lambda_0^\frac12 \bold{P}\\
\lambda_1^\frac12\bold{A}(\bold{m}^k)
\end{bmatrix} 
\bold{u}^{k+1} = 
\begin{bmatrix}
\lambda_0^\frac12[\bold{d}+\bold{d}^k]\\
\lambda_1^\frac12[\bold{b+b}^k] 
\end{bmatrix},
\end{equation}
where $\bold{d}^k$ and $\bold{b}^k$ are the components of the dual variable $\bold{s}^k$ associated with the observation-equation and wave-equation constraints, and are formed by the running sum of the data and source residuals in iteration (see equation~\ref{scale_main1}).
%
The closed-form expression of the reconstructed wavefield is given by
\begin{equation}  \label{solve_sub_uu} 
\bold{u}^{k+1} = \big[ \lambda_0 \bold{P}^T\bold{P}+ \lambda_1 \bold{A}(\bold{m}^k)^T \bold{A}(\bold{m}^k)\big]^{-1}\big[\lambda_0  \bold{P}^T [\bold{d}+\bold{d}^k]+\lambda_1\bold{A}(\bold{m}^k)^T [\bold{b+b}^k] \big].
\end{equation}
The reconstructed wavefield can be computed numerically with linear algebra methods (direct or iterative methods) suitable for sparse matrices as reviewed by \citet{vanLeeuwen_2016_PMP}.
\subsubsection{Solving for $\bold{m}$}
In order to solve the bound-constrained TV regularized nonlinear problem described by equation \ref{main4b0}, we first tackle the nonlinearity issue by considering the special structure of the Helmholtz operator $\bold{A}$ given in equation \ref{helmholtz}.
Using the following approximation
\begin{eqnarray}
\bold{A}(\bold{m})\bold{u}^{k+1} &=& \bold{\Delta u}^{k+1} + \omega^2  \bold{C} (\bold{m})\text{diag}(\bold{m})\bold{B} \bold{u}^{k+1}, \label{nonlinear_op20} \nonumber \\
&\approx & \bold{\Delta}\bold{u}^{k+1} + \underbrace{\omega^2  \bold{C} (\bold{m}^k)\text{diag}(\bold{B} \bold{u}^{k+1})}_{\bold{L}(\bold{u}^{k+1})} \bold{m}, \label{nonlinear_op21}
\end{eqnarray}
we linearise the operator $\bold{A}$ with respect to $\bold{m}$ by building 
the matrix $\bold{C}$ from $\bold{m}^k$ to manage potential nonlinear boundary conditions. Note that this linearization step is not necessary when the PML absorbing conditions are used. We also exploit the bilinearity of the wave equation to permute $\bold{B}\bold{u}$ and $\bold{m}$ in the operator $\bold{L}(\bold{u}^{k+1})$.
%
With these manipulations, the sub-problem associated with $\bold{m}$ is recast as 
\begin{equation} \label{new_m0}
\bold{m}^{k+1} ~= ~ \underset{\bold{m}\in \mathcal{C}}{\arg\min} ~~ 
 \|\bold{m}\|_\text{TV}+\frac{\lambda_1}{2}\| \bold{L}(\bold{u}^{k+1})\bold{m} -\bold{y}^{k}\|_2^2. 
\end{equation}
where 
\begin{equation} \label{y}
\bold{y}^{k} = \bold{b} + \bold{b}^{k} -\bold{\Delta u}^{k+1}.
\end{equation}
%
%
We solve the constrained optimization problem described by equation \ref{new_m0} with the split Bregman method \citep{Goldstein_2009_SBM} to decouple the TV minimization subproblem (first term in equation~\ref{new_m0}) from the $\ell_{2}$ subproblem  (second term in equation~\ref{new_m0}) and force the box constraint. We introduce the auxiliary variables $\bold{p}_0$, $\bold{p}_1$ and $\bold{p}_2$ to perform this splitting, in which these variables being related to $\bold{m}$ by means of a simple equality constraint.
The auxiliary variable $\bold{p}_0$ is used to enforce the box constraint, while $\bold{p}_1$ and $\bold{p}_2$ are the variables of the TV minimization problem. 
Let us define $(\bold{p}_1 \; \bold{p}_2)$ as a two-column matrix, $\|(\bold{p}_1 \; \bold{p}_2)\|=\sqrt{\bold{p}_1^2 + \bold{p}_2^2}$ be a vector which contains $\ell_2$ norm of each row of $(\bold{p}_1 \; \bold{p}_2)$, and
$\sum\|(\bold{p}_1 \; \bold{p}_2)\|$ be the mixed $\ell_{2,1}$ norm ($\ell_1$ norm of $\|(\bold{p}_1 \; \bold{p}_2)\|$), which promotes sparsity. Then
by defining the objective function $J$ as
\begin{equation} \label{cost}
J(\bold{m},\bold{p}) = \sum\|(\bold{p}_1 \; \bold{p}_2)\| + \frac{\lambda_1}{2}\|\bold{L}(\bold{u}^{k+1})\bold{m} -\bold{y}^{k} \|_2^2,
\end{equation} 
the unconstrained optimization problem described by equation \ref{new_m0} can be written in a split and constrained form as
\begin{subequations}
\label{Split00}
\begin{eqnarray} 
&& \underset{{\bold{m},\bold{p},\bold{p}_0 \in \mathcal{C}}}{\arg\min}  ~~~~~~ J(\bold{m},\bold{p}) \\
&&\text{subject to} ~~~~
\bold{p}=\bold{\nabla}\bold{m},
\end{eqnarray}
\end{subequations}
where 
\begin{equation}
\bold{\nabla}=\begin{bmatrix} \bold{I} \\ \bold{\nabla}_1 \\ \bold{\nabla}_2 \end{bmatrix} \in \mathbb{R}^{3N\times N}, \hspace{1cm}
\bold{p}=\begin{bmatrix} \bold{p}_0 \\ \bold{p}_1 \\ \bold{p}_2 \end{bmatrix}\in \mathbb{R}^{3N\times 1}, \nonumber
\end{equation}
and $\bold{I}$ is the identity matrix.

The scaled augmented Lagrangian function (Appendix A) for the problem defined by equation \ref{Split00} is
\begin{equation} \label{AL2}
\mathcal{L}_A(\bold{m},\bold{p},\bold{\bar{q}}) = J(\bold{m},\bold{p})  + \frac12\|\bold{\nabla}\bold{m} - \bold{p}-\bold{\bar{q}}^k\|_{\bold{\Gamma}}^2-\frac{1}{2}\|\bold{\bar{q}}^k\|_{\bold{\Gamma}}^2,
\end{equation}
where $\bold{\bar{q}}$ is the scaled Lagrangian multipliers and
\begin{equation} \label{Gamma}
\bold{\Gamma} =
 \left[
    \begin{array}{r@{}c|c|c@{}l}
  &    \begin{smallmatrix}
       \gamma_0 & & 0 \\
          &\ddots&\\
        0 & & \gamma_0 \rule[-1ex]{0pt}{2ex}
      \end{smallmatrix} 
  & \mbox{\huge0}  
   & \mbox{\huge0} & \rlap{\kern5mm$N$} \\
      \hline
  &    \mbox{\huge0} &  
       \begin{smallmatrix}\rule{0pt}{2ex}
        \gamma_1 & & 0 \\
          &\ddots&\\
        0 & & \gamma_1 \rule[-1ex]{0pt}{2ex}
      \end{smallmatrix}    & \mbox{\huge0} &  \rlap{\kern5mm $N,$}  \\
            \hline
  &    \mbox{\huge0} &  \mbox{\huge0} &
       \begin{smallmatrix}\rule{0pt}{2ex}
        \gamma_2 & & 0 \\
          &\ddots&\\
        0 & & \gamma_2
      \end{smallmatrix}    &  \rlap{\kern5mm $N$} 
    \end{array} 
\right]
\end{equation}
with the penalty parameters $\gamma_0$, $\gamma_1$, $\gamma_2>0$. 
Again, applying the method of multipliers to find the saddle point of the problem~\ref{AL2} gives
\begin{subequations} \label{main2s}
\begin{eqnarray}
\bold{m}^{k+1},\bold{p}^{k+1} &=&  \underset{{\bold{m},\bold{p},\bold{p}_0 \in \mathcal{C}}}{\arg\min}  ~ J(\bold{m},\bold{p})  + \frac12\|\bold{\nabla}\bold{m} - \bold{p}-\bold{q}^k\|_{\bold{\Gamma}}^2 \label{main2sa} \\
\bold{q}^{k+1} &=& \bold{q}^k + \bold{p}^{k+1}-\nabla \bold{m}^{k+1}.  \label{main2sd}
\end{eqnarray} 
 \end{subequations}
where the bar of $\bold{\bar{q}}^k$ is removed for simplicity.\\
Substituting the explicit expression of $J$, equation \ref{cost}, into equation \ref{main2sa} leads to the following PRSM iteration:  
\begin{subequations}
\label{main3}
\begin{eqnarray}
\bold{m}^{k+1} &= & \underset{\bold{m}}{\arg\min} ~~ 
 \frac{\lambda_1}{2}\|\bold{L}(\bold{u}^{k+1})\bold{m}  -\bold{y}^{k} \|_2^2 +\frac{1}{2}  \|\bold{\nabla m}-\bold{p}^k-\bold{q}^k\|_{\bold{\Gamma}}^2 \label{main3b} \\
 \bold{q}^{k+\frac12} &=& \bold{q}^k + \frac{1}{2}[\bold{p}^{k}-\nabla \bold{m}^{k+1}].  \label{main3d0} \\
 \bold{p}^{k+1} &= & \underset{\bold{p}_0\in \mathcal{C},(\bold{p}_1,\bold{p}_2)}{\arg\min} ~~ 
  \sum\|(\bold{p}_1 \; \bold{p}_2)\|+\frac{1}{2}  \|\nabla \bold{m}^{k+1} -\bold{p}-\bold{q}^{k+\frac12}\|_{\bold{\Gamma}}^2 \label{main3c} \\
\bold{q}^{k+1} &=& \bold{q}^{k+\frac12} + \frac{1}{2}[\bold{p}^{k+1}-\nabla \bold{m}^{k+1}].  \label{main3d}
\end{eqnarray}
\end{subequations} 
Note that the weight $\frac{1}{2}$ in equations \ref{main3d0} and \ref{main3d} has a similar role as $\alpha$ in equations \ref{main4b00} and \ref{main4b1}.
%
Now we come up with a linear inverse subproblem for $\bold{m}$, equation~\ref{main3b}. Accordingly, the update $\bold{m}^{k+1}$ is obtained by solving the following system of linear equations in a least-squares sense:
\begin{equation} \label{new_mm00}
\begin{bmatrix}
\lambda_1^{\frac12}\bold{L}(\bold{u}^{k+1}) \\
\bold{\Gamma}^\frac12\nabla 
\end{bmatrix}
\bold{m}^{k+1} = 
\begin{bmatrix}
  \lambda_1^{\frac12}[\bold{b}+\bold{b}^k -\bold{\Delta u}^{k+1}]  \\
 \bold{\Gamma}^\frac12[\bold{p}^k+\bold{q}^k]
  \end{bmatrix},
\end{equation}
where we have substituted $\bold{y}^{k}$ by its explicit expression, equation~\ref{y}. In equation~\ref{new_mm00}, the first line
describes the information carried out by the reconstructed wavefield to update $\bold{m}$ via the wave-equation rewriting, while the second line 
describes the action of the TV regularization and bound constraints on $\bold{m}$  via its linear relation with the auxiliary variable $\bold{p}$. \\
The closed-form expression of $\bold{m}$ is given by
\begin{equation} \label{new_mm000}
\bold{m}^{k+1} = 
\big[
\lambda_1\bold{L}(\bold{u}^{k+1})^T \bold{L}(\bold{u}^{k+1}) + \nabla^T\bold{\Gamma}\nabla 
\big]^{-1} 
\big[\lambda_1\bold{L}(\bold{u}^{k+1})^T [\bold{b}+\bold{b}^k-\bold{\Delta u}^{k+1}] + \nabla^T\bold{\Gamma}[\bold{p}^k+\bold{q}^k]\big].
\end{equation}
As for the linear system \ref{solve_sub_uu}, $\bold{m}$ can be computed numerically with any suitable sparse linear algebra method. \\
The sub-problem for $\bold{p}$, equation \ref{main3c}, is straightforward to solve.
The objective function is separable with respect to the variable  $\bold{p}_0$ and the variables  $\bold{p}_1$ and  $\bold{p}_2$ (i.e., the optimization can be performed for $\bold{p}_0$ and $\bold{p}_1$, $\bold{p}_2$ separately).
The variable  $\bold{p}_0$ is solution of the following linear inverse problem
 \begin{equation} \label{projp0}
 \bold{p}^{k+1}_0 = \underset{\bold{p}_0\in \mathcal{C}}{\arg\min} \frac{\gamma_0}{2}  \|\bold{m}^{k+1} -\bold{p}_0-\bold{q}_0^k\|_2^2.
 \end{equation}
%
The $i$th element of the solution, $\bold{p}^{k+1}_0(i)$, is the closest element of $\bold{m}^{k+1}(i)-\bold{q}^k_0(i)$ to the desired set $[\bold{m}_l(i), \bold{m}_u(i)]$. 
Therefore, 
 \begin{equation} \label{projp0_1}
 \bold{p}^{k+1}_0 = \text{proj}_{\mathcal{C}} (\bold{m}^{k+1} - \bold{q}^k_0),
 \end{equation}
where the projection operator is
 $\text{proj}_{\mathcal{C}} (\bullet) = \min(\max(\bullet,\bold{m}_{l}), \bold{m}_{u}) $.\\
%
The variables  $\bold{p}_1$ and  $\bold{p}_2$ are updated via the following proximity operator:
\begin{equation} \label{proxim}
\bold{p}_1^{k+1},  \bold{p}_2^{k+1} = \underset{(\bold{p}_1 \; \bold{p}_2)}{\arg\min} ~~ 
  \sum \|(\bold{p}_1 \; \bold{p}_2)\| +\frac{\gamma}{2} 
   \| (\bold{p}_1 \; \bold{p}_2) - (\bold{z}_1 \; \bold{z}_2)\|_2^2
\end{equation}
where $\gamma=\gamma_1=\gamma_2$, and $\bold{z}_i= \nabla_i\bold{m}^{k+1}-\bold{q}_i^k ~ i=1,2$.
Proximity operators are generalization of projection operators \citep{Combettes_2011_PRO}. Equation \ref{proxim} describes a separable optimization problem with respect to $\bold{p}_1$ and  $\bold{p}_2$. Furthermore,  $\bold{p}_1$ and  $\bold{p}_2$ have closed-form expressions \citep{Goldstein_2009_SBM}
\begin{equation}
\label{proj}
   \bold{p}_i^{k+1}
    = \text{prox}_{\gamma}(\bold{z}_i),~~i=1,2,
\end{equation} 
where
\begin{equation} \label{prox0}
\text{prox}_{\gamma}(\bold{z}_i) = \frac{\bold{z}_i}{\|(\bold{z}_1 \; \bold{z}_2)\|} \max(\|(\bold{z}_1 \; \bold{z}_2)\| - \gamma,0).
\end{equation}
It can be seen that, for a single vector, the proximity operator defined in equation \ref{prox0} reduces to soft thresholding.

Considering all the above-mentioned processes, a pseudocode for the BTV-regularized IR-WRI algorithm is summarized in Algorithm \ref{Alg2cont0}.
Note that the lines 6-11 of the algorithm correspond to one ADMM iteration of the model-parameter updating. These operations could be iterated in an inner loop to update the BTV-regularized model several times after the wavefield reconstruction at each outer iteration. However, we observed numerically that only a single iteration of the inner loop guarantees the most efficient convergence of the full algorithm. This property has been noticed by \citet{Goldstein_2009_SBM} and is discussed more extensively in the framework of IR-WRI by \citet{Aghamiry_2019_IWR}. The inefficiency of the inner iterations can be understood by the fact that the original nonlinear problem is solved with an alternating-direction strategy (managed by the outer loop). This implies that each subproblem is solved from potentially inaccurate passive variables, this inaccuracy preventing an efficient minimization of the objective during inner iterations.
Furthermore, in order to drive the algorithm, we assumed that the constraint is feasible. 
However, it has been shown that in the case of infeasible linear constraints the ADMM iteration can still produce approximate solutions
that are stable \citep{Frick_2011_MPA,Jiao_2016_ADM}.
\begin{algorithm}
\caption{
BTV regularized IR-WRI algorithm based on the PRS algorithm.}
\label{Alg2cont0}
\small
Initialize: set the RHS errors $k=0, \bold{s}^0 = \bold{0}$, $\bold{q}^0 = \bold{0}$\\
Input: $\bold{m}^0$ (initial model parameters)\\
\While {convergence criteria  not satisfied}{
 $\bold{u}^{k+1} \leftarrow  \text{update according to}~~ \ref{solve_sub_uu}$\\
 $\bold{s}^{k+\frac{1}{2}} \leftarrow  \bold{s}^{k} + 0.5[\bold{s}-\bold{F(m}^k)\bold{u}^{k+1}]$\\
 $\bold{m}^{k+1} \leftarrow  \text{update according to}~~ \ref{new_mm00}$\\
 $\bold{q}^{k+\frac12} \leftarrow  \bold{q}^k +0.5[\bold{p}^k-\bold{\nabla}\bold{m}^{k+1}]$\\
 $\bold{p}_0^{k+1} \leftarrow \text{proj}_{\mathcal{C}} (\bold{m}^{k+1} - \bold{q}_0^{k+\frac12})$\\
 $\bold{p}_1^{k+1} \leftarrow  \text{prox}_{\gamma} (\nabla_1\bold{m}^{k+1}-\bold{q}_1^{k+\frac12})$\\
 $\bold{p}_2^{k+1} \leftarrow  \text{prox}_{\gamma} (\nabla_2\bold{m}^{k+1}-\bold{q}_2^{k+\frac12})$\\
 $\bold{q}^{k+1} \leftarrow  \bold{q}^{k+\frac12} +0.5[\bold{p}^{k+1}-\bold{\nabla}\bold{m}^{k+1}]$\\
 $\bold{s}^{k+1} \leftarrow  \bold{s}^{k+\frac{1}{2}}+ 0.5[\bold{s}-\bold{F(m}^{k+1})\bold{u}^{k+1}]$\\
 $k \leftarrow k+1$ \\
}
\end{algorithm}
%
%

\section{Numerical examples}

\subsection{Experimental setup and parameter tuning}
We assess the performance of our BTV regularized IR-WRI against 2D mono-parameter synthetic examples. We start with a toy example built with a high-velocity inclusion model that is embedded in a background medium where velocity linearly increases with depth . To tackle more realistic applications, we proceed with two targets of the challenging 2004 BP salt model \citep{Billette_2004_BPB}. With the  BP salt case study, we seek to illustrate the potential of IR-WRI equipped with our BTV regularization to image salt bodies and sub-salt structures starting from crude initial models and realistic frequencies. \\
For all the numerical examples, forward modelling is performed with a 9-point stencil implemented with anti-lumped mass and PML absorbing boundary conditions \citep{Chen_2013_OFD}. In this setting the diagonal matrix $\bold{C}$ contains the damping PML coefficients and does not depend on $\bold{m}$.  With this setting, equation~\ref{nonlinear_op20} does not require any approximation for linearization. \\
%
We will compare the results of WRI and IR-WRI to highlight the improved convergence history of IR-WRI resulting from the iterative updating of the right-hand sides in the penalty function associated with the scaled-form augmented Lagrangian, equations~\ref{new_scale}-\ref{main4}. We assume that our IR-WRI algorithm, when this right-hand side updating is not activated, is representative of the WRI penalty method \citep{VanLeeuwen_2013_MLM}. For a fair comparison, we will use the same experimental setup (penalty parameters and stopping criterion of iteration) for the two methods. We also compare the WRI and IR-WRI results when they are obtained without any priors ($\gamma_0=\gamma=0$), with only bound constraints ($\gamma=0$) and with BTV  regularization ($\gamma_0 \ne 0$ and $\gamma \ne 0$), where it is reminded that $\gamma_0$ and $\gamma$ are the penalty parameters that control the weight of the bound constraints and TV regularization, respectively, in the objective function (see equation~\ref{Gamma}).  \\
We tune the different penalty parameters according to the following guideline. We start from the last subproblem of the splitting procedure and set the parameter $\gamma$, which controls the soft thresholding performed by the TV regularization, equation~\ref{proj}-\ref{prox0}. In this study, we find that $\gamma=2\% \max {\|(\bold{z}_1 \; \bold{z}_2)\|}$ was a good pragmatical value. This tuning can be refined according to prior knowledge of the geological structure, coming  from well log for example. In this study, we use the same weight for the bound constraints and the TV regularization: $\gamma_0=\gamma$. Once we set $\gamma$, we define $\lambda_1$ such that $\gamma / \lambda_1$ is a percentage of mean absolute value of the diagonal coefficients of $\bold{L}^T\bold{L}$ during the parameter estimation subproblem, equation \ref{new_mm000}. This percentage is set according to the weight that we want to assign to the TV regularization and bound constraints relatively to the wave equation constraint during the parameter estimation. Parameter $\lambda_1$ may be increased during iterations to reduce the weight of TV regularization and bound constraints near the convergence point. We found this adaptation useful when we start from very crude initial models. Finally, we set $\lambda_0$ such that  $\lambda=\lambda_1/\lambda_0$ is a small fraction of the highest eigenvalue $\xi$ of the normal operator $\bold{A(m})^{-T}\bold{P}^{T}\bold{P}\bold{A(m})^{-1}$ during the wavefield reconstruction subproblem, equation~\ref{solve_sub_uu}, according to the criterion proposed by \citet{vanLeeuwen_2016_PMP}. In all the numerical tests, we use  $\lambda=1\text{e-5} \xi$ and $\lambda=1\text{e-3} \xi$ for noiseless and noisy data, respectively. This tuning of $\lambda$ is indeed important because it controls the extension of the search space. 
A too high value of $\lambda$ reduces the weight of $\|\bold{Pu-d}\|_2^2$ during the wavefield reconstruction and makes IR-WRI to behave like a reduced approach. Conversely, using a small value for $\lambda$ fosters data fitting and expends the search space accordingly. However, a too small value can lead to a prohibitively high number of iterations of the augmented Lagrangian method before the wave equation constraint is fulfilled with a sufficient accuracy. Moreover, when data are contaminated by noise, a too small value for $\lambda$ will make the wavefield reconstruction to over-fit the data and drive WRI toward poor minimizer.
We always use $\lambda$ as a fixed percentage of $\xi$ in iterations for both WRI and IR-WRI. This does not prevent IR-WRI to converge towards accurate minimizers thanks to the iterative error correction performed by the Lagrange multiplier updating. The reader is referred to \citet{Aghamiry_2019_IWR} for a more thorough sensitivity analysis of IR-WRI to the penalty parameter $\lambda$.

\subsection{Inclusion model}
\begin{figure} 
\centering
   \includegraphics[width=0.5\textwidth]{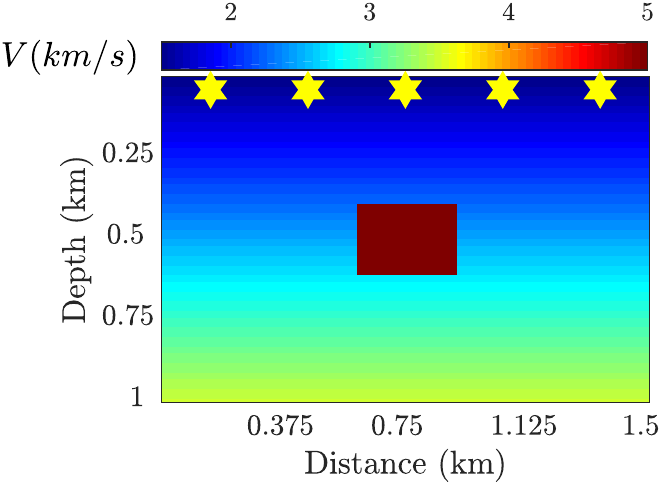} 
\caption{Box-shape anomaly example. True velocity model. The yellow stars show the source positions.}
\label{fig_Box2_true}
\end{figure}

\begin{figure} 
\centering
   \begin{subfigure}[b]{1\textwidth}
   \includegraphics[width=1\textwidth]{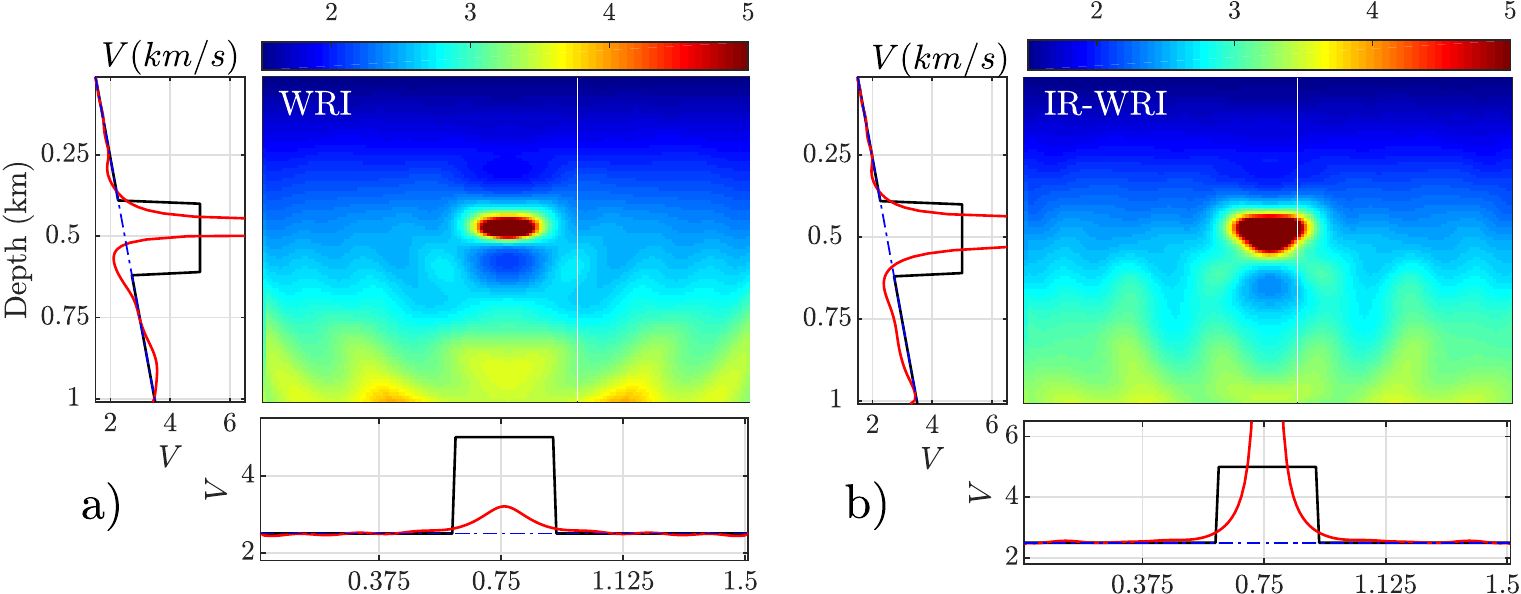} 
\end{subfigure}
\begin{subfigure}[b]{1\textwidth}
   \includegraphics[width=\textwidth]{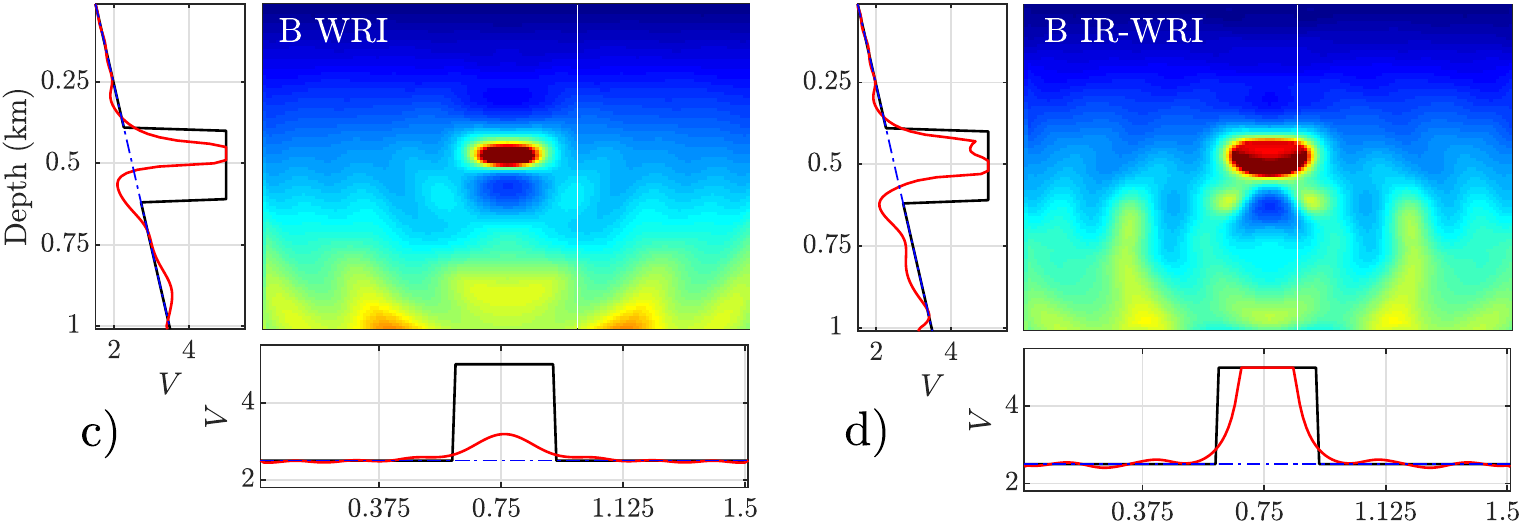}
\end{subfigure}
\begin{subfigure}[b]{1\textwidth}
   \includegraphics[width=\textwidth]{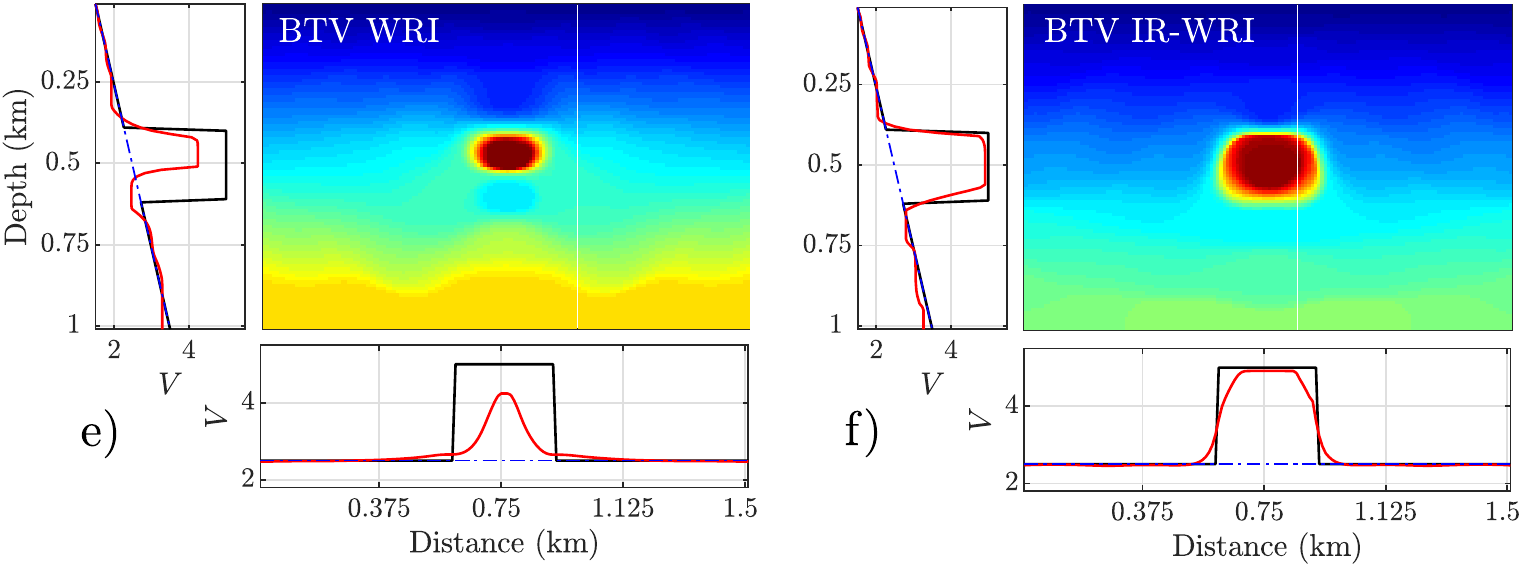}
\end{subfigure}
\caption{Box-shape anomaly example with the true background velocity gradient as initial model. (a) WRI. (b) IR-WRI. (c-d) Bound constrained WRI (c) and IR-WRI (d) models. (e-f) BTV regularized WRI (e) and IR-WRI (f) models. Horizontal and vertical profiles across the center of the inclusion from the true (black), initial (dash blue) and reconstructed (red) models are shown below and on the left-hand side of the models.}
\label{fig_Box2_new}
\end{figure}

The subsurface model contains a sharp box-shape anomaly of sides 0.2 $\times$ 0.3~km with a velocity ($V_P$) of 5~km/s. It is embedded in a smooth background model where $V_P$ increases linearly with depth from 1.5 to 3.5~km/s (Figure~\ref{fig_Box2_true}). The model is 1.5~km long and 1~km deep, and is discretized with a 10~m grid interval. The regular surface acquisition consists of five sources (as depicted with yellow stars in Figure~\ref{fig_Box2_true}) and 65 receivers deployed on the surface. The source signature is a Ricker wavelet with a 5~Hz dominant frequency. We start the inversion from the true background model and invert simultaneously three frequency components (2.5, 5 and 7 Hz) with noiseless data. 
This frequency bandwidth has been selected to cover a significant band of vertical wavenumbers in the waveform-inversion sensitivity kernels, considering the limited aperture illumination provided by the surface acquisition. Moreover, a realistic starting frequency of 2.5Hz allows us to assess the resilience to cycle skipping of IR-WRI.
A maximum number of iterations set to 70 is used as stopping criterion for all of the tests shown in Figure~\ref{fig_Box2_new}. When bound constraints are used, the bounds $\bold{m}_{l}$ and $\bold{m}_{u}$, equation~\ref{projp0_1}, is set to the true minimum and maximum square slownesses, respectively. 

\begin{figure} 
\centering
   \begin{subfigure}[b]{1\textwidth}
   \includegraphics[width=1\textwidth]{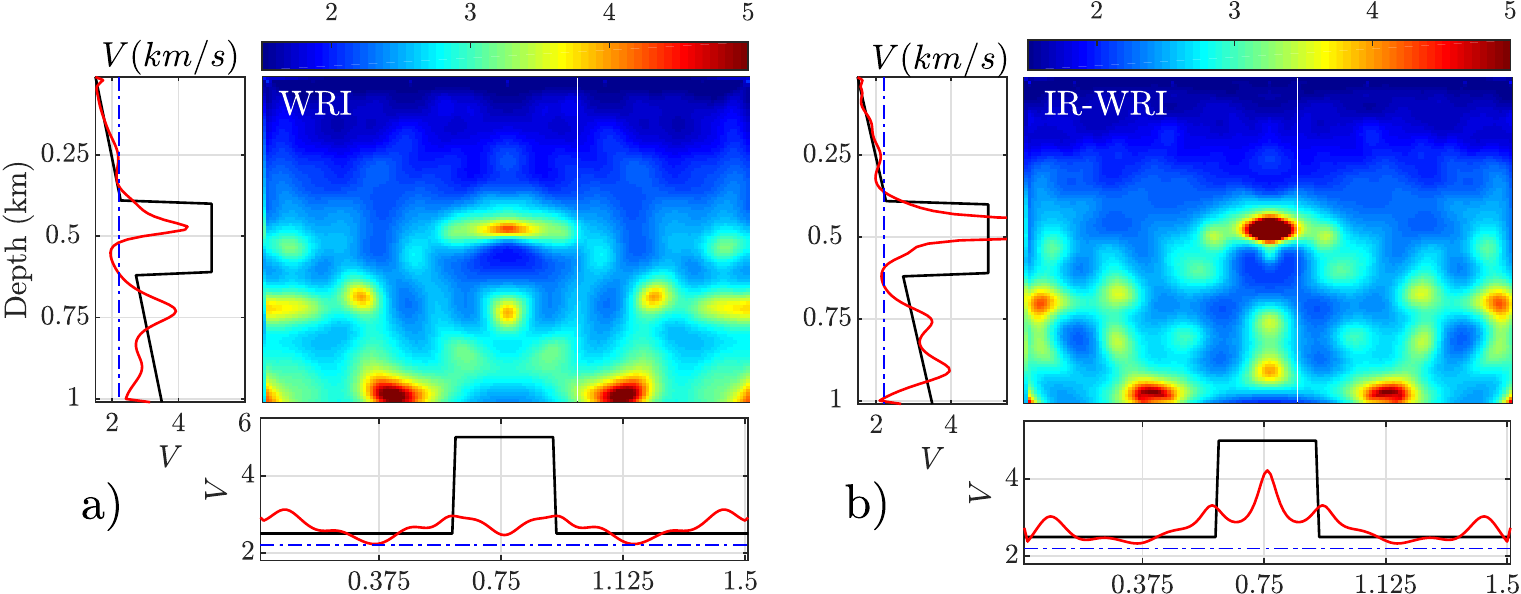} 
\end{subfigure}
\begin{subfigure}[b]{1\textwidth}
   \includegraphics[width=\textwidth]{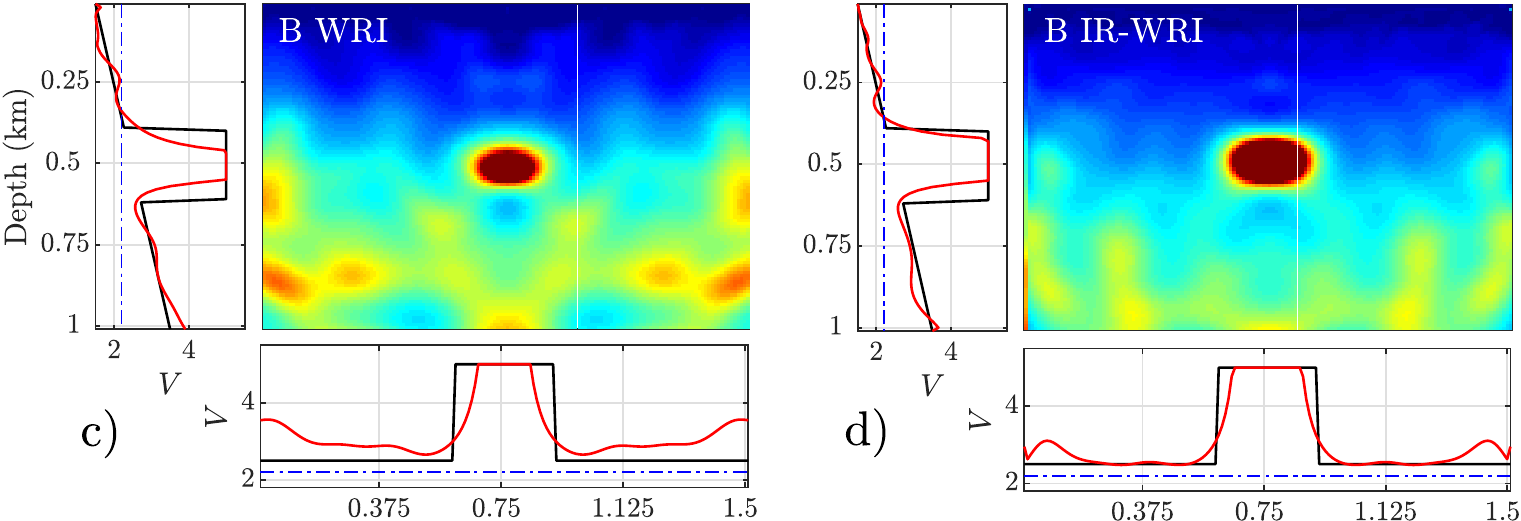}
\end{subfigure}
\begin{subfigure}[b]{1\textwidth}
   \includegraphics[width=\textwidth]{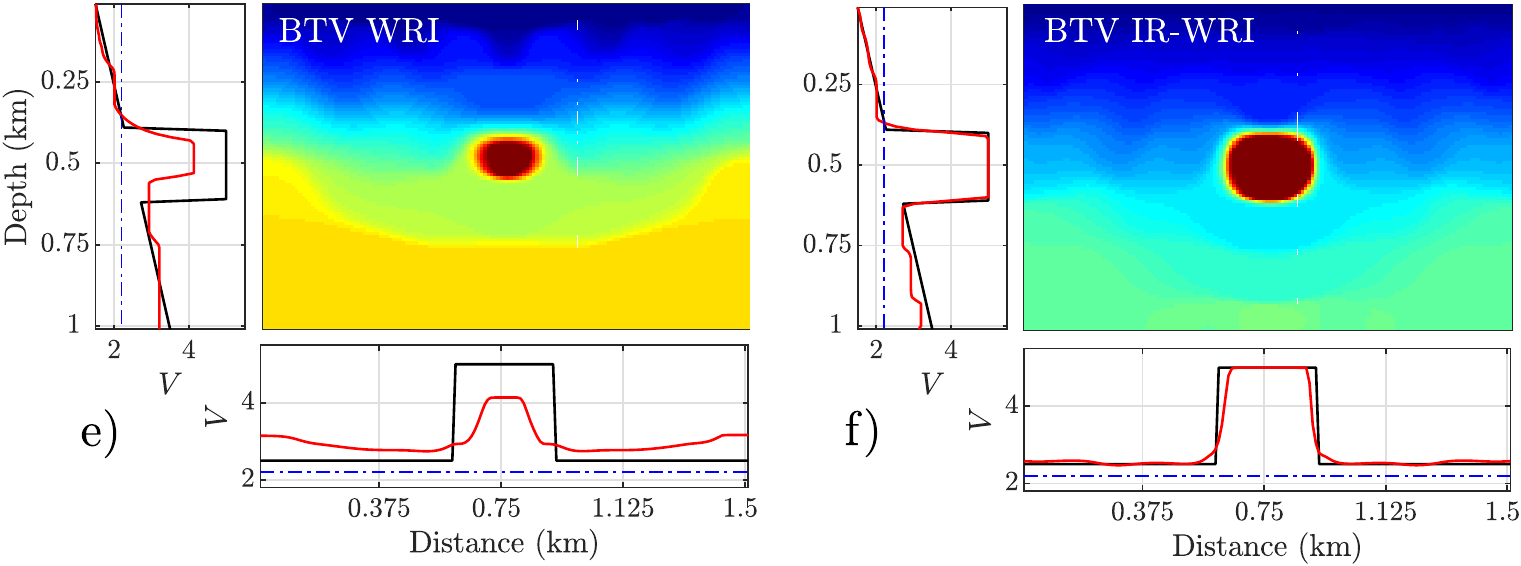}
\end{subfigure}
\caption{Box-shape anomaly example. Same as Figure~\ref{fig_Box2_new} for a homogeneous initial velocity model ($V_P$=2.2km/s).}
\label{fig_Box_new}
\end{figure}

\begin{figure} 
\centering
   \begin{subfigure}[b]{0.7\textwidth}
   \includegraphics[width=\textwidth]{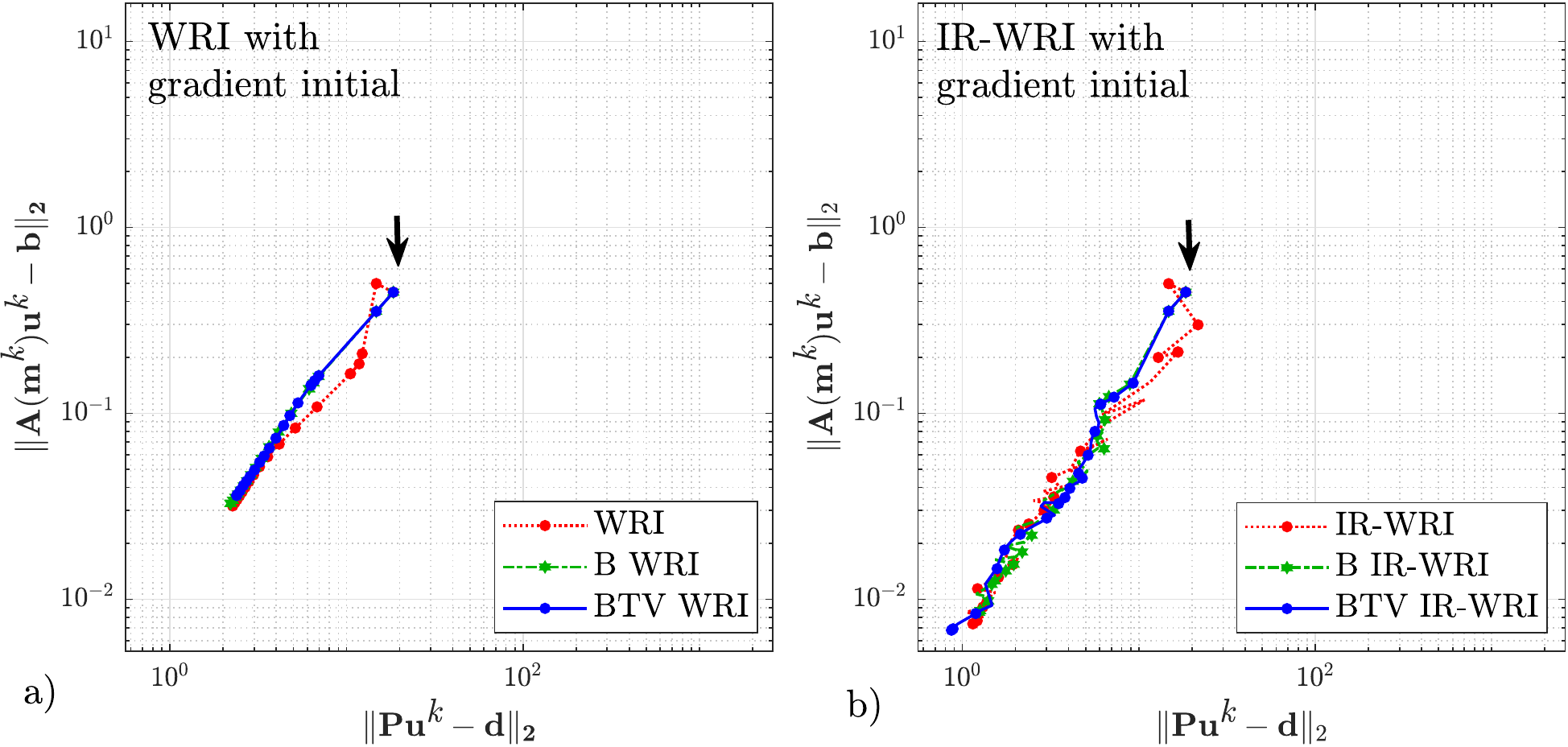} 
\end{subfigure}
\begin{subfigure}[b]{0.7\textwidth}
   \includegraphics[width=\textwidth]{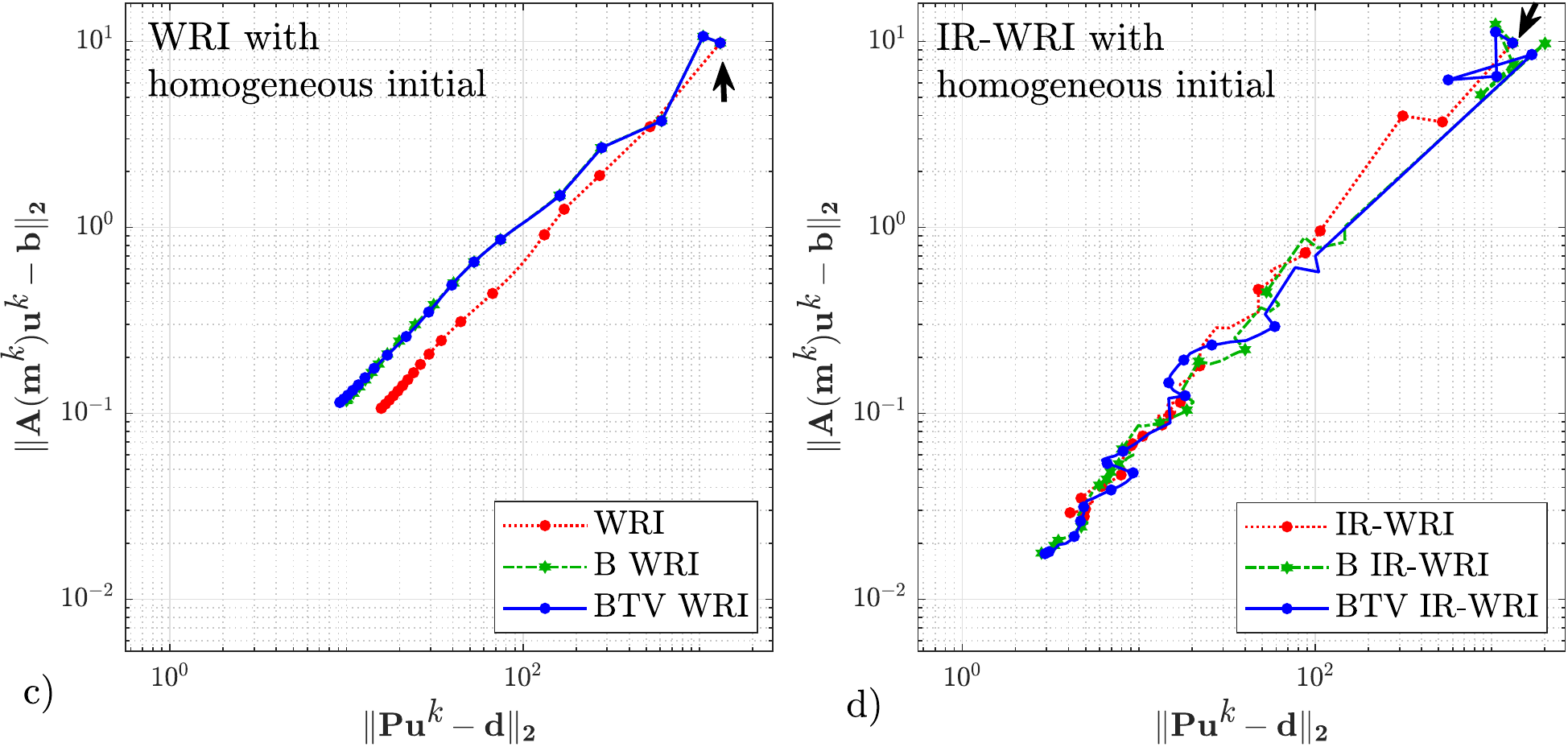}
\end{subfigure}
\caption{Box-shape anomaly example. Convergence history  in the ($\|\bold{Pu}^k-\bold{d}\|_2 - \|\bold{A(m}^k)\bold{u}^k-\bold{b}\|_2$) plane of WRI and IR-WRI without priors, with bound constraints and with BTV regularizations. (a,c) WRI. (b,d) IR-WRI. Initial model is (a-b) the true velocity-gradient background model, and (c-d) the homogeneous velocity model.  All the panels are plotted with the same horizontal and vertical logarithmic scale. The black arrow points the starting point.}
\label{fig_Box_res}
\end{figure}

\begin{figure} 
\centering
   \begin{subfigure}[b]{0.7\textwidth}
   \includegraphics[width=1\textwidth]{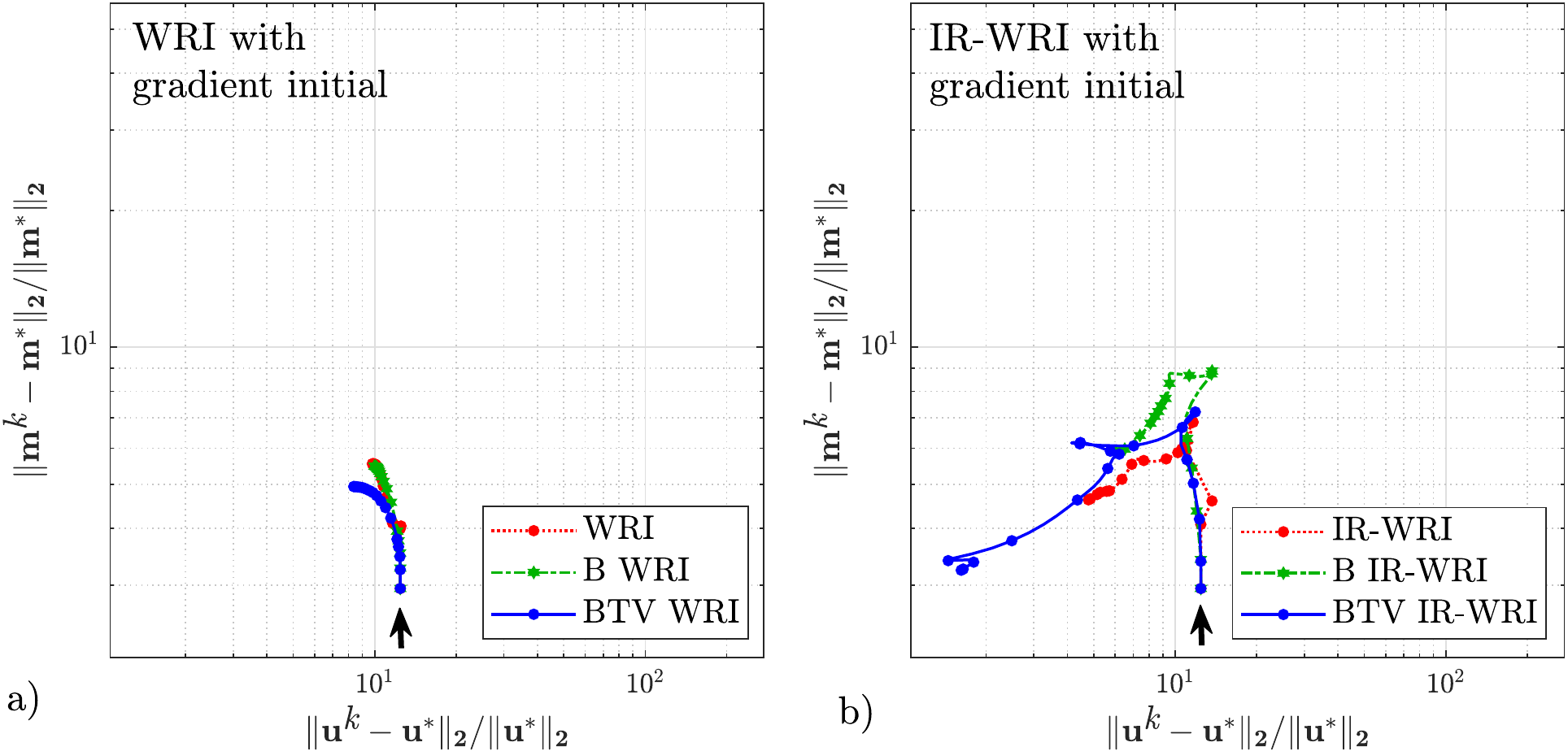} 
\end{subfigure}
\begin{subfigure}[b]{0.7\textwidth}
   \includegraphics[width=\textwidth]{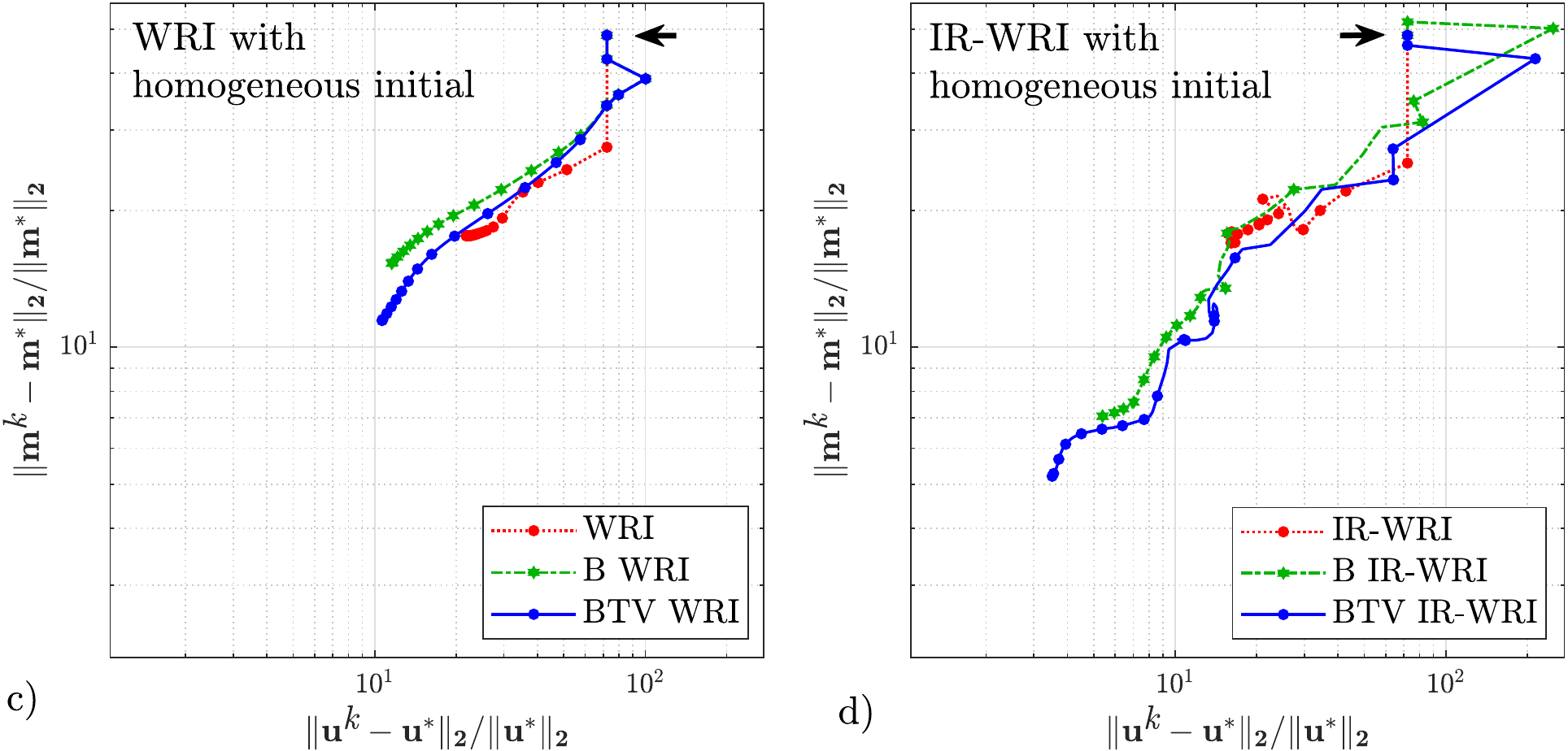}
\end{subfigure}
\caption{Box-shape anomaly example. Convergence history  in the ($\|\bold{u}^k-\bold{u}^*\|_2/\|\bold{u}^*\|_2-\|\bold{m}^k-\bold{m}^*\|_2/\|\bold{m}^*\|_2$) plane of WRI and IR-WRI without priors, with bound constraints and with BTV regularizations and for the two initial models. (a,c) WRI. (b,d) IR-WRI. Initial model is (a-b) the true velocity-gradient background model, and (c-d) the homogeneous velocity model. Note the increase of $\|\bold{m}^k-\bold{m}^*\|_2/\|\bold{m}^*\|_2$ over iterations (see text for explanations). All the panels are plotted with the same horizontal and vertical logarithmic scale. The black arrow indicates the starting point.}
\label{fig_Box_mse}
\end{figure}

\begin{figure} 
\centering
   \includegraphics[width=0.7\textwidth]{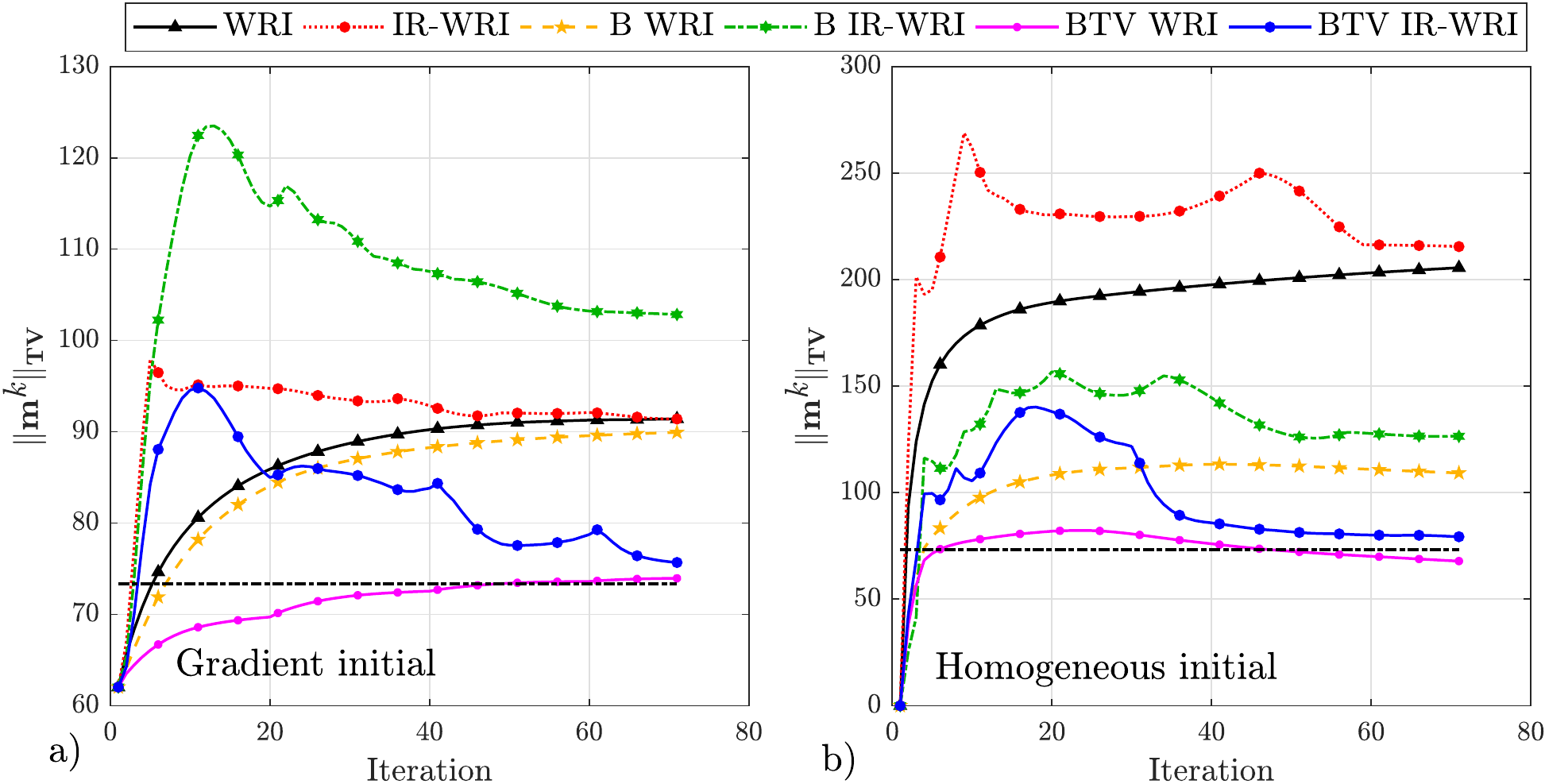} 
\caption{Box-shape anomaly example. TV norm history over iterations of  bound constrained, BTV regularized and ordinary  WRI and IR-WRI for (a) velocity-gradient and (b) homogeneous initial model. The TV norm of true model is plotted with dashed-black line.}
\label{fig_Box_tv}
\end{figure}

We first compare WRI and IR-WRI results when bound constraints and TV regularization are not activated, i.e. $\gamma_0=\gamma=0$ (Figure~\ref{fig_Box2_new}(a-b)). WRI and IR-WRI reconstruct only the top of the anomaly with strongly overestimated velocities. Then, we add bound constraints in WRI and IR-WRI, using $\gamma_0/\lambda_1=0.01 \zeta$, where $\zeta$ is the mean absolute value of the diagonal coefficients of $\bold{L}^T\bold{L}$.  
The bound-constrained WRI and IR-WRI only reconstruct the top of the anomaly as in Figure~\ref{fig_Box2_new}(a-b). However, the inclusion velocities are now well controlled by the bound constraints (Figure~\ref{fig_Box2_new}(c-d)). These first two tests show that IR-WRI reconstructs better the shape of the anomaly than WRI with however more significant artifacts on both sides of the anomaly. These artifacts may result from the deficit of horizontal-wavenumber illumination provided by the sparse limited-offset surface acquisition and by multi-scattering pollutions. \\
%
Then, we apply BTV regularization with $\gamma / \lambda_1=\gamma_0 / \lambda_1=0.01 \zeta$ (Figure~\ref{fig_Box2_new}(e-f)).  
Since the initial model matches the true velocity-gradient background model, we use a small value of $\gamma / \lambda_1$ (i.e., a high value of $\lambda_1$) and keep it constant in iterations to preserve the smooth components of the subsurface model. 
We show that the BTV regularized WRI still fails to reconstruct the full anomaly (Figure~\ref{fig_Box2_new}e). In contrast, BTV regularized IR-WRI keeps on improving the reconstruction of the anomaly in depth, while efficiently mitigating the oscillating artifacts (Figure~\ref{fig_Box2_new}f). \\
%
A vertical profile across the reconstructed anomaly also highlights some limitations of the BTV regularization (Figure~\ref{fig_Box2_new}f):  below the anomaly, the BTV regularization superimposes staircase artifacts  on the velocity gradient, consistently with the piecewise constant assumption underlying TV regularization.  \\
%
%
%
%
To emphasize the resilience to cycle skipping of the BTV regularized WRI and IR-WRI, we repeat this toy example using a 2.2~km/s homogeneous velocity model as initial model (Figure~\ref{fig_Box_new}).
We perform a first test without any prior. Compared to the previous test, we just stabilize the inversion by adding a small damping term  (=$0.01 \zeta$) to $\bold{L}^T\bold{L}$.
Compared to Figure~\ref{fig_Box2_new}(a-b), the artifacts have a much stronger imprint due to the inaccuracy of the starting model (Figure~\ref{fig_Box_new}(a-b)).
Then, we move to bound-constrained and TV regularized tests (Figure~\ref{fig_Box_new}(c-f)). To decrease the above-mentioned artifacts, we assign a high initial weight to the BTV regularization ($\gamma / \lambda_1=\gamma_0 / \lambda_1=\zeta$) and decrease it by a factor 2 every 10 iterations until it reach a minimal value set to $0.01 \zeta$ (i.e., the constant value previously used ).
In accordance with the former test,  bound constraints alone are not sufficient to reconstruct the bottom part of the anomaly and cancel out the oscillating artifacts (Figure~\ref{fig_Box_new}(c-d)). In contrast, BTV IR-WRI achieves these two goals, although it leaves a significant staircase footprint below the anomaly  (Figure~\ref{fig_Box_new}f). Compared to Figure ~\ref{fig_Box2_new}f, the edges of the anomaly are better reconstructed at the expense of the background velocity-gradient model. This results because more aggressive TV regularization was used during the early iterations of this test, allowing for a better reconstruction of the blocky components of the medium, while injecting undesired staircase footprint on its smooth components.
As for the former test, IR-WRI clearly outperforms WRI due to the more efficient solution refinement procedure resulting from the right-hand side updating. \\
%
%
%
We also show the joint evolution in iterations of the observation-equation and wave-equation errors (Figure~\ref{fig_Box_res}) and the wavefield and subsurface model errors (Figure~\ref{fig_Box_mse}), when the initial model is the true velocity-gradient background model and the homogeneous model.
IR-WRI fits the data and wave equation better than WRI after 70 iterations with both initial models because the right-hand side updating embedded in IR-WRI cancels out more efficiently the data and source residuals in iterations and, hence better refines the solution accordingly. Moreover, BTV regularization in IR-WRI further improves the data and wave equation fit for both initial models because it reduces more efficiently the oscillating artifacts in the reconstructed velocity model and better reconstructs the edges of the anomaly (Figure~\ref{fig_Box2_new}f). 
As above mentioned, the oscillating artifacts may result from the deficit of horizontal wavenumber illumination generated by the limited-offset surface acquisition. In this framework, the prior contained in the BTV regularization efficiently narrows the null space of the inversion.
The more complex zigzag path followed by IR-WRI relative to WRI in the ($\|\bold{Pu}^k-\bold{d}\|_2 - \|\bold{A(m}^k)\bold{u}^k-\bold{b}\|_2$) plane highlights how the joint updating of the data and source by their associated residuals dynamically balances the weight of the two objective functions in iterations. This more complex convergence history of IR-WRI, which has been already noticed in \citet{Aghamiry_2019_IWR}, suggests that the right-hand side updating perform a self-adaptive weighting of the two competing objective functions driven by the relative reduction of the data and source residuals in iterations. 
This zigzag convergence trend translates also into more complex path in the $(\|\bold{u}^k-\bold{u}^*\|_2/\|\bold{u}^*\|_2-\|\bold{m}^k-\bold{m}^*\|_2/\|\bold{m}^*\|_2)$ plane ($\bold{m}^*$ and $\bold{u}^*$ denote the true model and wavefield, respectively), which illustrates how the solution refinement is pushed toward the wavefield reconstruction or the velocity model estimation according to the self-adaptive weighting of the observation-equation and wave-equation objective functions (Figure~\ref{fig_Box_mse}). Note that the relative model errors increase in the case of the initial velocity-gradient model (Figure~\ref{fig_Box_mse}a,b). This results because the smooth background model is degraded by the oscillating artifacts and the staircase footprint during the sharp inclusion reconstruction. Indeed, this degradation of the smooth components has a much higher weight in the $\ell_{2}$ misfit function than the more accurate reconstruction of the blocky components. However, this increase is much more moderate in IR-WRI (Figure~\ref{fig_Box_mse}b) than in WRI (Figure~\ref{fig_Box_mse}a).
Finally, we plot the TV norm of the reconstructed models in iterations for the different tests (Figure \ref{fig_Box_tv}). As expected, the models reconstructed with bound constraints and TV regularization match better the TV of the true model. The BTV WRI model matches slightly better the TV of the true model than the BTV IR-WRI one, in particular when the initial model is the true background model. Indeed, this does not reflect that WRI better reconstructs the anomaly than IR-WRI. Instead, it reflects the slower convergence of WRI relatively to IR-WRI which contributes to keep the background model smooth (namely, which a TV close to 0 and equal to that of the true model).
%

%

%
\subsection{2004 BP salt model - central target}
\begin{figure} 
\centering
   \begin{subfigure}[b]{1\textwidth}
   \includegraphics[width=1\textwidth]{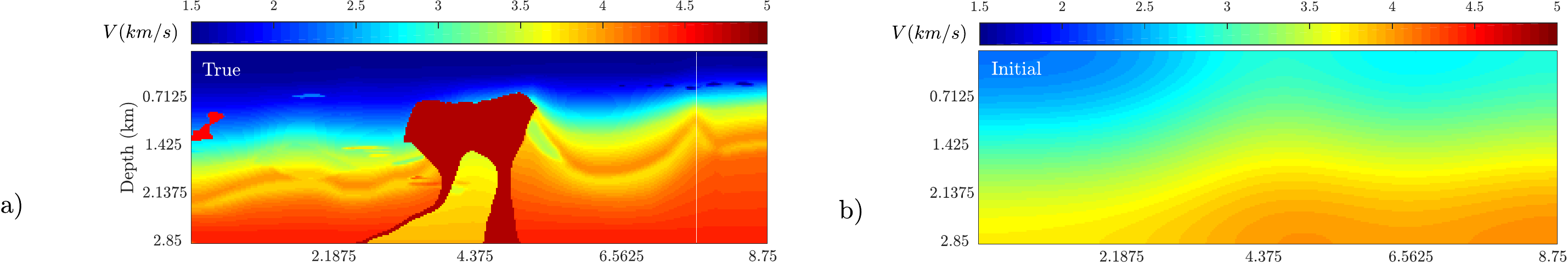} 
\end{subfigure}
\begin{subfigure}[b]{1\textwidth}
   \includegraphics[width=\textwidth]{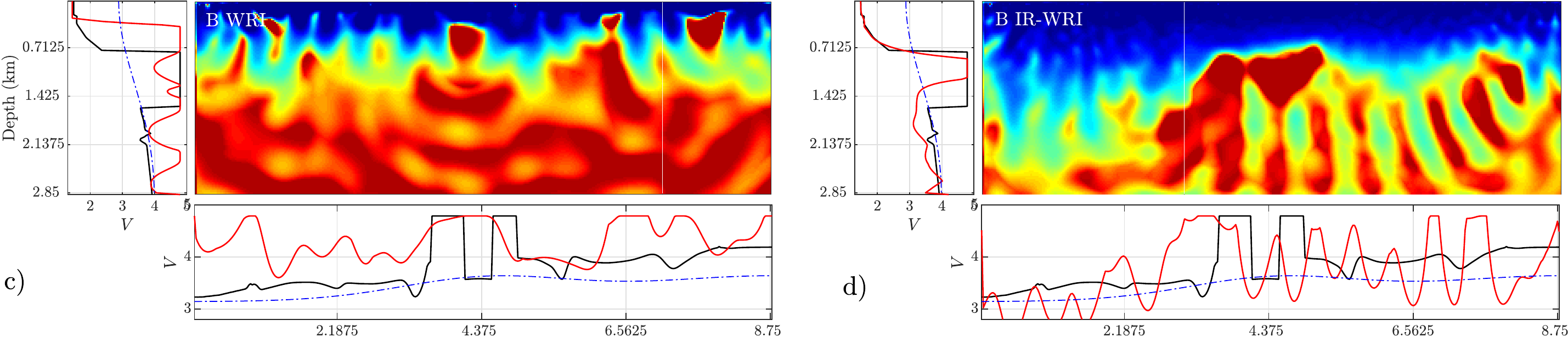}
\end{subfigure}
\begin{subfigure}[b]{1\textwidth}
   \includegraphics[width=1\textwidth]{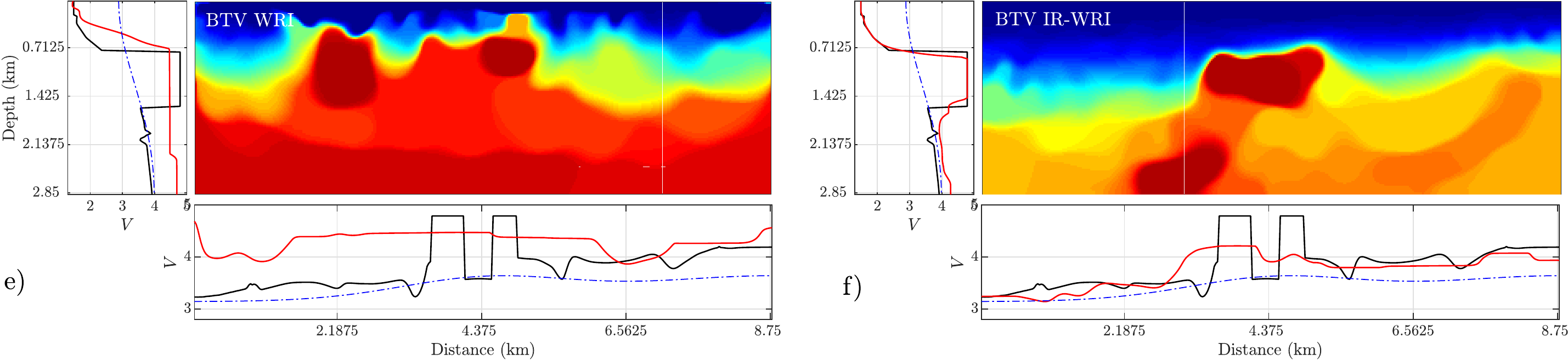} 
\end{subfigure}
\caption{2004 BP salt model - Central target. 3-Hz frequency. (a) True model. (b) Starting model. (c-d) Bound constrained WRI(c) and IR-WRI(d) models. (e-f) BTV regularized WRI(e) and IR-WRI(f) models. Horizontal and vertical profiles at 4.15~km distance and 1.65~km depth from the true (black), initial (dash blue) and reconstructed (red) models are shown below and on the left-hand side of the models.}
\label{fig_Bp_part_new}
\end{figure}

\begin{figure}
\begin{center}
\includegraphics[scale=0.5]{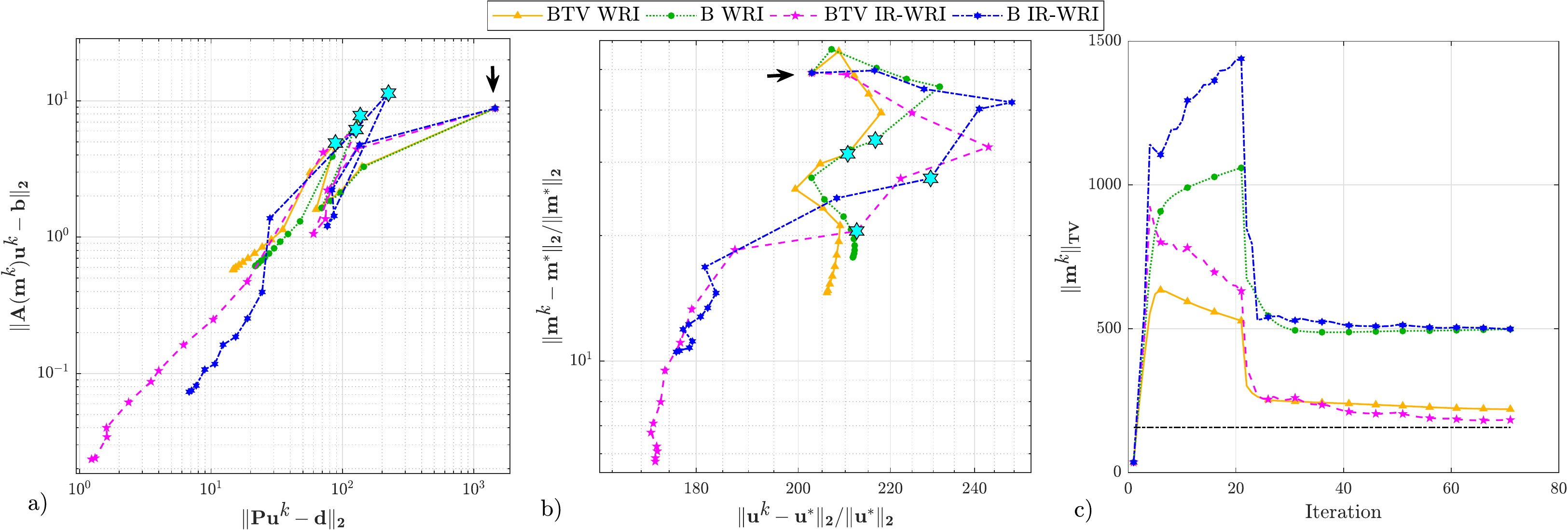}
\caption{2004 BP salt model - Central target. 3-Hz frequency. Convergence history (a) in the ($\|\bold{Pu}^k-\bold{d}\|_2 - \|\bold{A(m}^k)\bold{u}^k-\bold{b}\|_2$) plane (b) in the ($\|\bold{u}^k-\bold{u}^*\|_2/\|\bold{u}^*\|_2-\|\bold{m}^k-\bold{m}^*\|_2/\|\bold{m}^*\|_2$) plane. (c) evaluation of TV norm $\|\bold{m}_TV\|$ over iterations. Note again the more complex convergence history of IR-WRI compared to WRI due to the self-adaptive weighting of the data-fitting and wave-equation objectives performed by right-hand side updating. The iteration 22 are located with cyan stars in (a) and (b).}
\label{fig_Bp_part_residuals_new}
\end{center}
\end{figure}

We now consider a more realistic application with a first target of the challenging 2004 BP salt model. The  2004 BP salt model is representative of the geology of the deep offshore Gulf of Mexico and mainly consists of a simple background with a complex rugose multi-valued salt body, sub-salt slow velocity anomalies related to over-pressure zones and a fast velocity anomaly to the right of the salt body \citep{Billette_2004_BPB}. The first selected target corresponds to the central part of the 2004 BP salt model characterized by a  deeply rooted salt body (Figure \ref{fig_Bp_part_new}a). \\
The subsurface model is 8.8~km wide and 2.9~km deep, and is discretized with a 25~m grid interval. We used 50 sources spaced 175~m apart on the top side of the model. The source signature is a Ricker wavelet with a 10~Hz dominant frequency. A line of receivers with a 50~m spacing are deployed at the surface leading to a stationary-receiver acquisition.  \\
We start the inversion from a smoothed version of true velocity model where the imprint of the salt body and other structures were cancelled out (Figure~\ref{fig_Bp_part_new}b) and invert the 3-Hz frequency with noiseless data. \\
We compare the results of WRI and IR-WRI with bound constraints and BTV regularization. To highlight the specific role of bound constraints, we activate them after 21 iterations. Since we start from a rough initial model, we set $\gamma_0/ \lambda_1 = \gamma / \lambda_1= \zeta$ and decrease them during iterations in a manner similar to the box-shape anomaly test with homogeneous starting model. Also, we add a damping term to $\bold{L}^T\bold{L}$ with a weight equal to $0.01 \zeta$ to further stabilize the inversion. We stop inversion after 70 iterations. The estimated models are shown in Figure \ref{fig_Bp_part_new} together with horizontal profiles at 1.65~km depth and vertical profiles at 4.15~km distance extracted from the true, initial and reconstructed models. 
As for the inclusion test, WRI fails to reconstruct the salt body and the subsalt structure because the data and source residuals are not re-injected in the right-hand sides of the penalty function at each iteration as in equation~\ref{main30}, leading to a stagnant convergence of the inversion (Figure~\ref{fig_Bp_part_new}c,e). When IR-WRI is applied with bound constraints alone, the reconstructed model is affected by noise with a periodic horizontal pattern. This noise likely results from the monochromatic nature of the inversion, multi-scattering within the salt body and limited illumination of the horizontal wavenumbers of the salt body leading to wraparound (Figure~\ref{fig_Bp_part_new}d). The BTV regularized IR-WRI mitigates efficiently this noise without degrading the resolution of the salt body and the sub-salt structures (Figure~\ref{fig_Bp_part_new}f). \\
Figure \ref{fig_Bp_part_residuals_new}(a) shows the joint evolution in iterations of the data misfit and the wave-equation error. As for the inclusion test, note the zigzag path followed by the IR-WRI objective functions over iterations. Also, the joint evolution of wavefield and subsurface model errors and the evolution of TV norm over iteration are shown in Figure \ref{fig_Bp_part_residuals_new}(b,c).  The TV norm evolution emphasizes how the bound constraints fasten the convergence of TV-regularized IR-WRI after iteration 20.  

%
%
\begin{figure} 
\centering
\begin{subfigure}[b]{1\textwidth}
   \includegraphics[width=\textwidth]{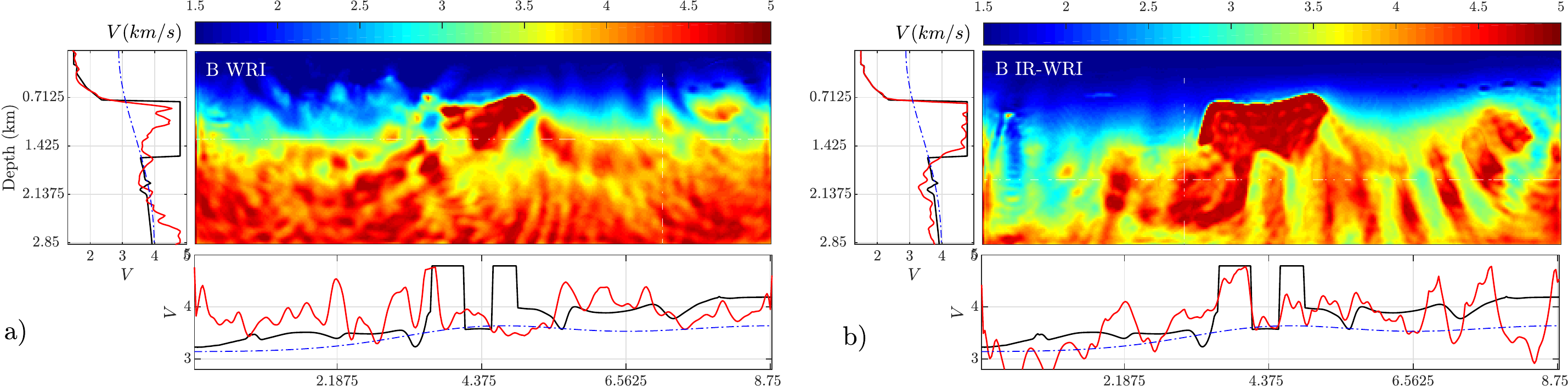}
\end{subfigure}
\begin{subfigure}[b]{1\textwidth}
   \includegraphics[width=1\textwidth]{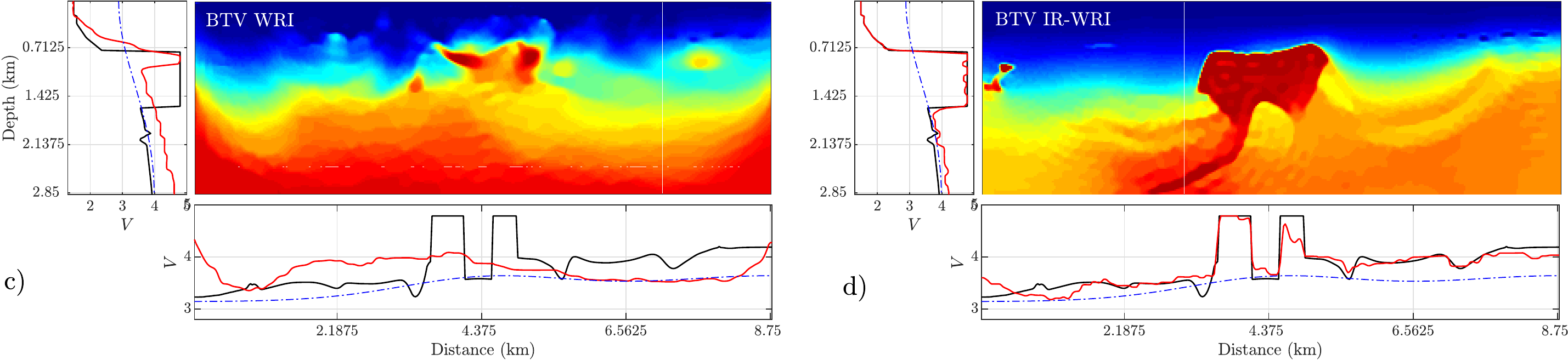} 
\end{subfigure}
\caption{2004 BP salt case study. Central target. (a) Final bound constrained WRI with Figure \ref{fig_Bp_part_new}(c) as initial model.(b) Final bound constrained IR-WRI with Figure \ref{fig_Bp_part_new}(d) as initial model. (c) Final BTV regularized WRI with Figure \ref{fig_Bp_part_new}(e) as initial model. (d) Final BTV regularized IR-WRI with Figure \ref{fig_Bp_part_new}(f) as initial model. }
\label{fig_Bp_final}
\end{figure}

We continue the inversion at higher frequencies using the final models of the 3~Hz inversion as initial model (Figure~\ref{fig_Bp_part_new}c-f). We used small batches of two frequencies with one frequency overlapping between two consecutive batches, moving from the low frequencies to the higher ones according to a classical frequency continuation strategy. We set $\gamma_0 / \lambda_1=\gamma / \lambda_1=0.01 \zeta$ and remove the damping of $\bold{L}^T\bold{L}$. The starting and final frequencies are 3.5~Hz and 12~Hz and the sampling interval in one batch is 0.5~Hz. The algorithm performs at most 15 iterations per frequency batch and the number of iterations that have been performed is 170. The inversion results are shown in Figure~\ref{fig_Bp_final}. WRI with bound constraints and BTV (Figure~\ref{fig_Bp_final}(a,c)) fails to converge toward satisfactory results, while BTV-regularized IR-WRI converges to accurate velocity model, although a significant imprint of the TV regularization is shown (Figure~\ref{fig_Bp_final}d). When IR-WRI is performed with only bound constraints, the oscillating artifacts are not cancelled out (Figure~\ref{fig_Bp_final}b). This highlights the role of TV regularization in reconstructing blocky structures and removing wraparound artifacts.

%
%
\subsection{2004 BP salt model - Left target}
\begin{figure}
\begin{center}
\includegraphics[scale=0.7]{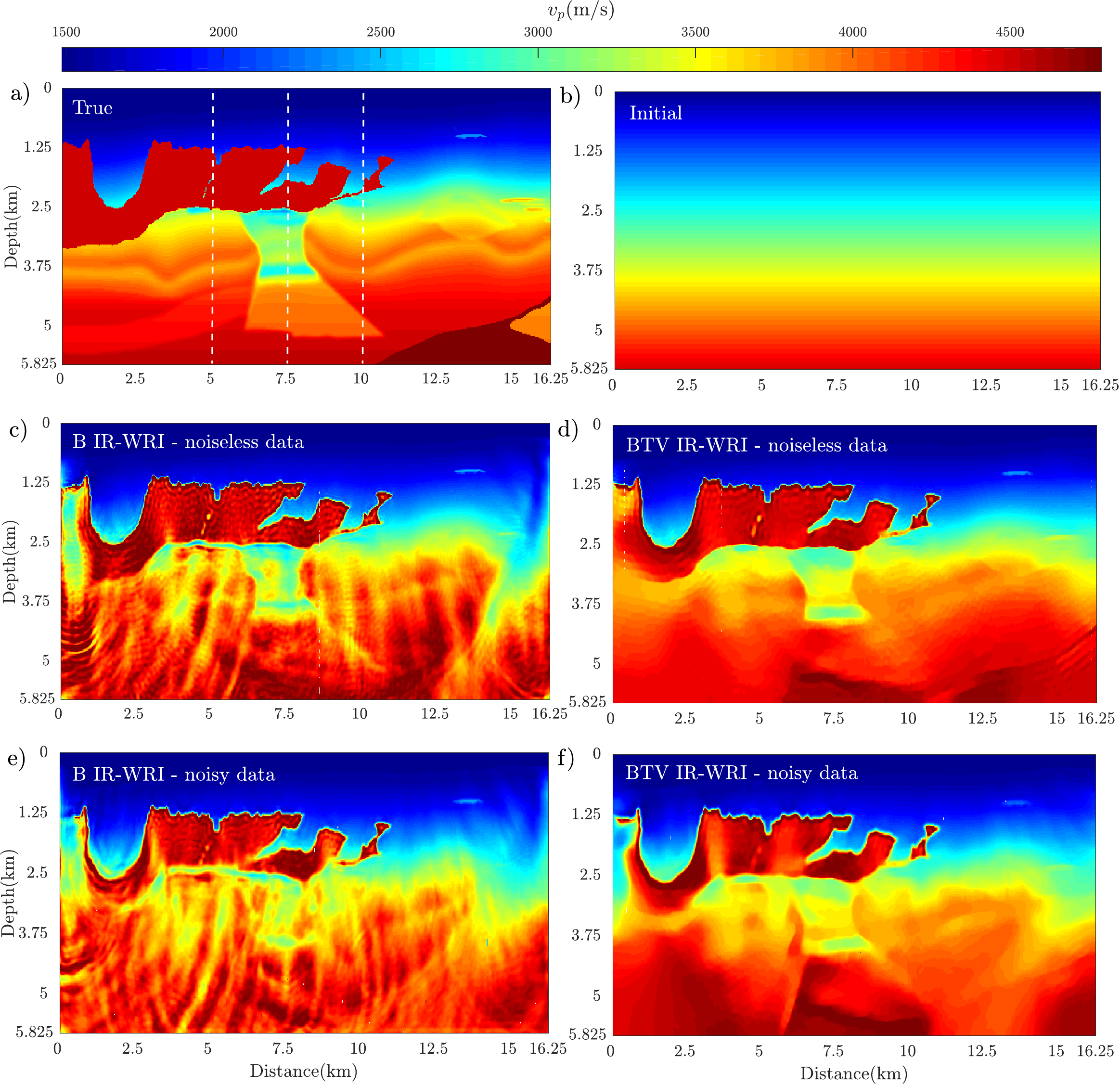}
\caption{2004 BP salt case study. Left target. (a) True BP model. The vertical dashed lines indicate the location of vertical logs of Figure~\ref{fig_Bp_logs}. (b) Initial velocity model. (c-d) Final IR-WRI velocity models obtained with bound constraints (c) and with BTV regularization (d) for noiseless data. (e-f) Same as (c-d) for noisy data for a SNR of 10~db.}
\label{fig_Bp}
\end{center}
\end{figure}

We now consider a second target of the 2004 BP salt model located on the left side of the model (Figure \ref{fig_Bp}a). This target was  previously used by among others \citet{Metivier_2016_TOF}, \citet{Brandsberg-Dahl_2017_FMU}, \citet{Esser_2018_TVR} and \citet{Kalita_2018_IBC} for FWI applications.  The subsurface model is 16250~m wide and 5825~m deep, and is discretized with a 25~m grid interval. We used 108 sources spaced 150~m apart on the top side of the model. The source signature is a Ricker wavelet with a 10~Hz dominant frequency. A line of receivers with a 25~m spacing are deployed at the surface leading to a stationary-receiver acquisition. 
We perform IR-WRI with bound constraints alone and with BTV regularization, for noiseless and noisy data. 
%

\begin{figure}
\begin{center}
\includegraphics[scale=0.7]{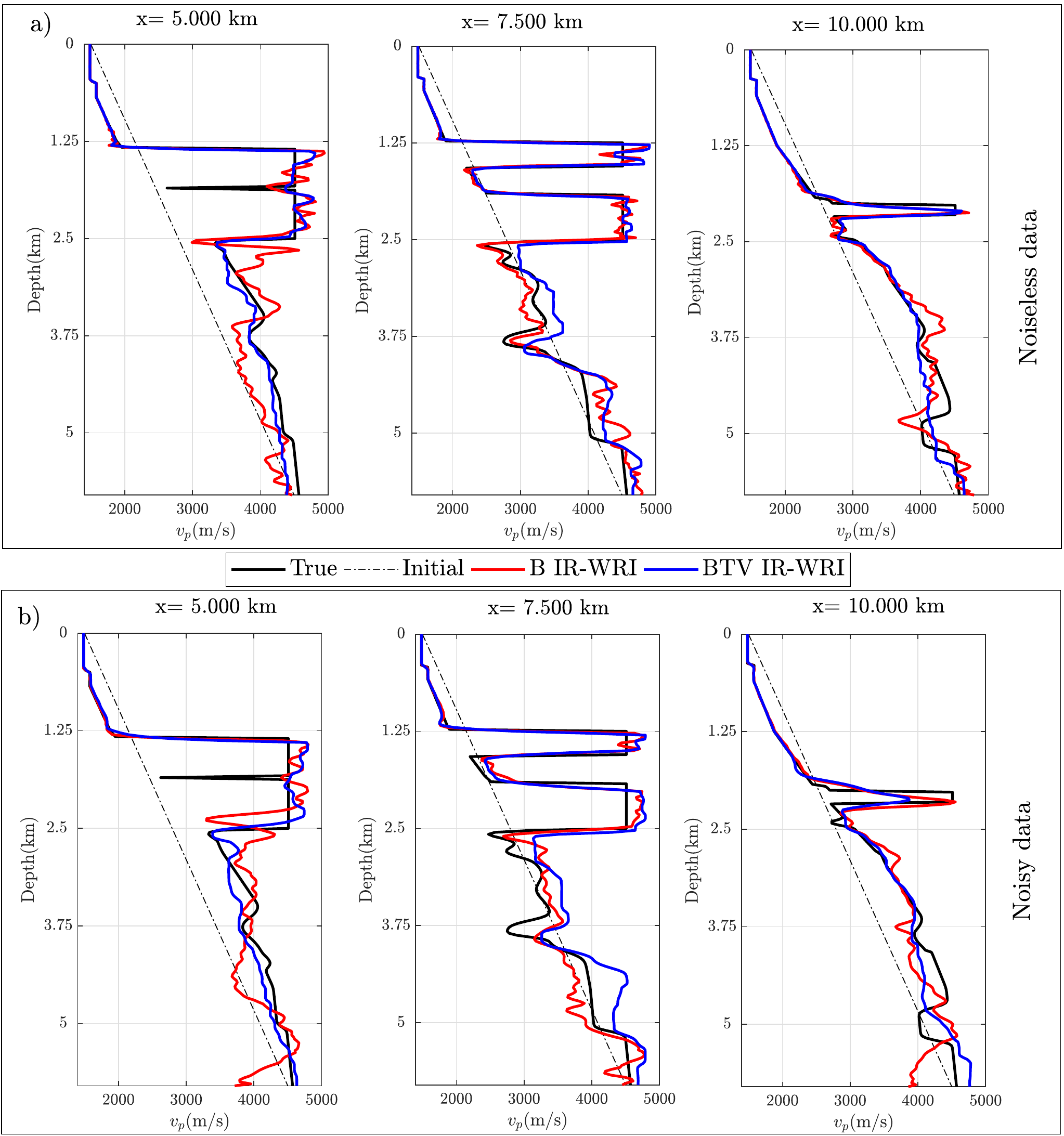}
\caption{2004 BP salt case study - Left target. Direct comparison between the true velocity model (black), the initial model (dashed line) and the IR-WRI models obtained with bound constraints (red) and with  BTV regularization (blue) along three logs at x=4.5~km, 7~km, 10~km (vertical dashed lines in Figure~\ref{fig_Bp}a)  from left to right. (a) Noiseless data. (b) Noisy data.}
\label{fig_Bp_logs}
\end{center}
\end{figure}

We used a crude laterally-homogeneous velocity-gradient model as initial model (Figure \ref{fig_Bp}b). 
We used small batches of two frequencies with one frequency overlap between two consecutive batches, moving from the low frequencies to the higher ones according to a classical frequency continuation strategy. The starting and final frequencies are 3~Hz and 13~Hz and the sampling interval in one batch is 0.5~Hz. The stopping criterion of iteration for each batch is given by
\begin{equation} \label{Stop}
k_{max}=15~~~~\text{or}~~~(\| \bold{A(m}^{k+1})\bold{u}^{k+1}-\bold{b}\|_2 \leq \varepsilon_b  \hspace{1.2 cm} \text{and}  \hspace{1.2 cm} \|\bold{Pu}^{k+1}-\bold{d}\|_2 \leq \varepsilon_d),
\end{equation} 
where $k_{max}$ denotes the maximum iteration count, $\varepsilon_b$=1e-3, and $\varepsilon_d$=1e-5 for noiseless data and $\varepsilon_b$=1e-3 , $\varepsilon_d$= noise level of batch for noisy data. 
We perform three paths through the frequency batches to improve the IR-WRI results, using the final model of one path as the initial model of the next one (these cycles can be viewed as outer iterations of IR-WRI). The starting frequency of the second and third path is 6~Hz and 8.5~Hz, respectively. 
The IR-WRI models inferred from noiseless data with bound constraints and with BTV regularization are shown in Figure~\ref{fig_Bp}(c-d). The number of iterations that have been performed with bound constraints and with BTV regularization are 441 and 340, respectively. Direct comparison between the true model, the starting model and the IR-WRI models along three vertical logs cross-cutting the salt body at 5~km, 7.5~km and 10~km distance (vertical dashed lines in Figure~\ref{fig_Bp}a) are shown in Figure~\ref{fig_Bp_logs}a. 
The results show the resilience of IR-WRI to cycle skipping with a pretty accurate reconstruction of the salt body and sub-salt structures (Figure~\ref{fig_Bp}(c-d)). 
However, the model obtained with bound constraints alone shows high-frequency noise in the salt body and below (Figure~\ref{fig_Bp}c). This noise can result from Gibbs phenomenon resulting from the frequency decimation and artifacts resulting from multi scattering. The BTV regularization efficiently removes these artifacts except those resulting from truncation of the acquisition near the left end of the model at 5-km depth (Figure~\ref{fig_Bp}d). The inversion captures reasonably well the subsalt structures, including the low-velocity over-pressure zone at (x,z)=(7.5km,4km) as well as the smooth velocity variations.  \\
When noisy data are used (Figure~\ref{fig_Bp}(e-f) and \ref{fig_Bp_logs}b), the number of iterations that have been performed with bound constraints alone and with BTV regularization is 263 and 254, respectively. 
As for the noiseless case, a direct comparison between the true model, the starting model and the two IR-WRI models along three vertical logs cross-cutting the salt body at 5~km, 7.5~km and 10~km distance is shown in Figure~\ref{fig_Bp_logs}b. The artifacts have now a more significant imprint when only bound constraints are used (Figure~\ref{fig_Bp}e). Accordingly, we apply a more aggressive TV regularization to obtain the results shown
in  Figure~\ref{fig_Bp}f. The artifacts have been efficiently removed with however a more obvious imprint of the piecewise-constant approximation underlying BTV regularization.  This blocky pattern is clearly visible in deep part of the vertical profiles of Figure~\ref{fig_Bp_logs}b where velocity gradients have been replaced by stack of constant-velocity layers.

%

\section{Discussion}
We have implemented bound constraints and TV regularization in the wavefield reconstruction inversion (WRI) method \citep{VanLeeuwen_2013_MLM} which has been recently improved by \citet{Aghamiry_2019_IWR} in the framework of the alternating-direction method of multiplier (ADMM), leading to the iteratively-refined WRI (IR-WRI).
To do this, we formulate IR-WRI as a TV minimization problem subject to observation-equation, wave-equation and bound constraints. We use a scaled form of the method of multiplier to recast 
the augmented Lagrangian function as  a quadratic penalty function in which the scaled Lagrange multipliers, here the running sum of the data and source residuals, are added to the right-hand sides (data and source) at each iteration.
At each iteration of the workflow, we perform a first ADMM step to break down the wavefield reconstruction and the parameter estimation into a sequence of two subproblems.
The wavefield reconstruction subproblem is a linear inverse problem which can be solved with direct or iterative methods for sparse linear systems \citep{Duff_2017_DMS,Saad_2003_IMS}.
Once one iteration of the wavefield reconstruction has been performed, we tackle the parameter estimation subproblem which involves the mixed $\ell_{1,2}$ TV norm of the model, the $\ell_{2}$ wave-equation objective and the bound constraints. We apply a second ADMM step to decompose this multi-objective optimization problem into a sequence of simpler subproblems. A first interesting property is the bilinearity of the wave equation constraint with respect to the wavefield and the parameter, which makes the $\ell_{2}$ wave-equation objective quadratic. A second key point is to introduce auxiliary variables, which are linearly related to the subsurface parameters and describes the TV of the model and the bound constraints. They allows for the de-coupling between the $\ell_{1}$ and $\ell_{2}$ components of the penalty function following the split Bregman method proposed by \citet{Goldstein_2009_SBM}. After this de-coupling, the subsurface parameter are first updated 
by minimizing the source residuals through the resolution of a sparse BTV-regularized Gauss-Newton system, before updating 
the auxiliary variables with proximal operators. This cycle is iterated until convergence. 

There are some key differences between our implementation of BTV-regularized WRI and previous ones based upon projected-gradient method \citep{Peters_2016_CVP} and primal dual hybrid gradient (PDHG) methods \citep{Esser_2018_TVR,Yong_2018_TVR} that we review below.
In the projected-gradient method, \citet{Peters_2016_CVP} first update model parameters with FWI or WRI  and then project the updated model into the intersection of the TV and bound constraint with Dykstra projection algorithm \citep{Boyle_1986_MFP}, an alternating projection onto some constraints until satisfying all of them. Therefore, their workflow is subdivided in two different parts, that are the model update  followed by the projection onto intersection of all constraints \citep[][ their equation 7]{Peters_2016_CVP}.    
Unlike \citet{Peters_2016_CVP} method which rely on independent update and projection steps, all the ingredients of BTV-regularized IR-WRI (i.e., wavefield reconstruction, parameter estimation, TV regularization and bound constraints) are consistently integrated in the theoretical framework of ADMM optimization \citep[the readers can also refer to ][]{Maharramov_2015_TVM}. \\
%
%
\citet{Esser_2018_TVR,Yong_2018_TVR} implement TV regularization and bound constraints in the reduced variable projection WRI of \citet{vanLeeuwen_2016_PMP} with PDHG. PDHG is a method to solve constrained optimization problems by alternating gradient descent (for primal variable) and gradient ascent (for dual variable), which can be interpreted as linearized ADMM \citep{Goldstein_2015_APD}. PDHG can be helpful if the least squares minimizations embedded in ADMM are difficult to solve efficiently \citep{Goldstein_2015_APD}. This is not really the case in WRI, which mainly requires to solve two sparse linear systems for $\bold{u}$ and $\bold{m}$. Meanwhile, the selection of step size in PDHG that guarantees fast convergence, or even convergence at all, is not intuitive at all and can make PDHG impractical. This issue prompted \citet{Goldstein_2015_APD} to develop step size tuning rules which contribute to make PDHG self-adaptive as illustrated recently by \citet{Yong_2018_TVR} in the frame of WRI.  Beyond PDHG method, \citet{Esser_2018_TVR,Yong_2018_TVR} minimize data and source residuals with a penalty method, which lacks the convergence property of the augmented Lagrangian method promoted in IR-WRI. 
Moreover, they implement TV regularization as a hard constraint in the parameter-estimation subproblem through a Lagrangian function, whose saddle point is estimated with PDHG.  
In contrast, we implement the TV regularization as a soft constraint after introducing the auxiliary variable $\bold{p}$ and solve the regularized sub-problem for $\bold{m}$ with the split Bregman method (or equivalently ADMM), equation~\ref{AL2}. This gives us the necessary flexibility to implement aggressive regularization during early iterations and relax it progressively, when very crude initial models are used. With more accurate initial models, the augmented Lagrangian embedded in the split Bregman method, equation~\ref{AL2}, allows for constant penalty parameter to be used. These penalty parameters are used as step lengths in the augmented Lagrangian method, hence leading to a self-adaptive TV regularization implementation rid of tedious TV-norm ball continuation strategies and/or adaptive step lengths \citep{Esser_2018_TVR,Yong_2018_TVR}. 
Third, while \citet[][ equation 18]{Esser_2018_TVR} and \citet[][ equation 41]{Yong_2018_TVR} implement bound constraints as hard constraints by projection at each iteration of the TV-regularized perturbation model onto the  feasible set defined by the bound constraints, we implement bound constraints consistently with TV regularization in the framework of the method of multiplier, equations~\ref{new_mm00}-\ref{proxim}.  Finally, our approach does not require prior guess of  TV-norm ball because we implement TV regularization as a minimization problem rather than as a constraint. That being said, \citet{Esser_2018_TVR} show promising results with the $\tau$ continuation approach and asymmetric TV norm on the BP salt model starting from a very crude initial model. It will be interesting in future work to assess the benefit of asymmetric TV norm in BTV IR-WRI. \\
Indeed, implementation of BTV-regularized IR-WRI requires to setup different parameters, which have been discussed at the beginning of the section {\it{Experimental setup and parameter tuning}}. The penalty parameter which controls the relative weight between the wave equation objective and the data fitting objective during the wavefield-reconstruction subproblem has been discussed in length in \citet{Aghamiry_2019_IWR}. \citet{Aghamiry_2019_IWR} have concluded that a fixed penalty parameter can be used during iterations because the accuracy of the minimizer in the method of multiplier is controlled both by the penalty parameter and the accuracy of the updated multiplier \citep[][ Theorem 17.6]{Nocedal_2006_NOO}. In this study, we suggest good pragmatical values of this penalty parameter as a percentage of the maximum eigenvalue of the augmented wave-equation normal operator for noiseless and noisy data.
Optimal values may be refined by trial and error to prevent noisy data over-fitting and keep the iteration count within reasonable limits.
The TV and bound parameters should be easily determined from well logs or a priori geological knowledge. Finally, the relative weight between the TV regularization and the wave-equation constraint, which is controlled by $\lambda_1$, needs also to be estimated during the parameter-estimation subproblem. As above mentioned, this penalty parameter is kept fixed during iterations or is progressively decreased to relax the TV regularization and bound constraints near the convergence point, depending of the accuracy of the initial model. The reader is also referred to \citet[][ section 2.2]{Goldstein_2009_SBM} who discuss the sensitivity  of the split Bregman method (an optimization method similar to ADMM \citep{Esser_2009_ALA}) to the penalty parameter for $\ell_1$-regularized problems. \\
One drawback of the BTV regularization is related to the piecewise constant approximation underlying TV regularization, which tends to superimpose some blocky patterns on the smooth part of the subsurface. To overcome this issue, ongoing work seeks to optimally combine Tikhonov and TV regularization in IR-WRI \citep{Gholami_2013_BCT,Aghamiry_2018_HTT}. \\
Other perspective works involve extension to multiparameter reconstruction, 3D geometries and application to real data to further assess the potential and limits of IR-WRI. For 3D applications, the wavefield reconstruction in the frequency domain requires to solve a large-scale linear algebra problem, equation~\ref{solve_sub_uu}. \citet{Operto_2015_ETF}, \citet{Amestoy_2016_FFF}, and \citet{Operto_2018_MFF} have shown the computational efficiency of 3D frequency-domain FWI based on sparse direct solvers in the 3.5-10 Hz frequency band for dense stationary-recording ocean-bottom cable acquisitions.  \citet{Mary_2017_PHD} showed that block low-rank multifrontal solver allows one to tackle numerical problems involving up to 100 million unknowns. The symmetry of the normal operator, equation~\ref{solve_sub_uu}, should balance the computational overhead resulting from the higher number of non-zero coefficients relative to the impedance matrix ${\bf{A}}$. Note that ${\bf{A}}$ was processed as an unsymmetric matrix in the above-mentioned references due to discretization issues, although more suitable discretizations which preserve the symmetry of the impedance matrix may be considered in the future \citep{Pratt_2018_FWI}. Alternatively, domain decomposition methods suitable for Helmholtz problems can be interfaced with hybrid direct-iterative solvers to perform wavefield reconstruction in bigger computational domains \citep{Dolean_2015_IDD}. These approaches may be suitable for stationary-receiver acquisition involving a more limited number of reciprocal sources sich as ocean bottom seismometer acquisitions. \\
Let's add that BTV IR-WRI should be easily implemented in the time-domain formulation of WRI recast as an extended-source waveform inversion \citep{Wang_2016_FIR,Huang_2018_VSE}.
%

\section{Conclusion}
We have presented a new method to implement TV regularization and bound constraints in frequency-domain FWI based on wavefield reconstruction (WRI). In a previous study, we have reformulated WRI in the framework of the alternating-direction method of multiplier (ADMM), leading to the iteratively-refined wavefield-reconstruction inversion method (IR-WRI). We have shown how the augmented Lagrangian embedded in ADMM makes IR-WRI weakly sensitive to a wide range of penalty parameter thanks to the Lagrange multiplier updating. Using a small value of this penalty parameter efficiently extends the search space during early iterations to foster data fitting without preventing the wave-equation constraint to be fulfilled at the convergence point with a preset prescribed error.
IR-WRI performs a first ADMM step to alternate wavefield reconstruction and  subsurface parameter estimation as in the original WRI method. When BTV regularization is used, we perform a second ADMM step to decompose the BTV-regularized parameter estimation sub-problem into a sequence of two simpler subproblems through the introduction of auxiliary variables. This variable splitting allows for the de-coupling between the $\ell_{1}$ and the $\ell_{2}$ components of the penalty function according to the so-called split Bregman method. An interesting property exploited by WRI methods is the bilinearity of the wave equation constraint with respect to the wavefield and the parameter, which makes the $\ell_{2}$ wave-equation objective of the second sub-problem quadratic.
Our implementation of BTV-regularization in IR-WRI with ADMM (or, equivalently split Bregman) provides a versatile framework to cascade constraints and regularization of different nature and is reasonably easy-to-tune due to the limited sensitivity of the augmented Lagrangian method to the choice of the penalty parameters. For challenging subsurface targets with contrasted structures such as salt bodies, we have shown that our BTV-regularized WRI shows a high resilience to cycle skipping and noise and efficiently mitigates high-frequency artifacts associated with incomplete illumination and multi scattering without detriment to the resolution of the imaging. 

\section{Acknowledgements}
This study was partially funded by the SEISCOPE consortium (\textit{http://seiscope2.osug.fr}), sponsored by AkerBP, CGG, CHEVRON, EXXON-MOBIL, JGI, PETROBRAS, SCHLUMBERGER, SHELL, SINOPEC, EQUINOR and TOTAL. This study was granted access to the HPC resources of SIGAMM infrastructure (http://crimson.oca.eu), hosted by Observatoire de la Cote d’Azur and which is supported by the Provence-Alpes C\^ote d'Azur region, and the HPC resources of CINES/IDRIS under the allocation 046091 made by GENCI.

\appendix
\section{Scaled form of augmented Lagrangian}
Let's start with the following constrained problem
\begin{subequations}
\label{app1}
\begin{eqnarray}
&& \min_{\bold{x}} ~~~~~~~~~\quad \|P(\bold{x})\|_2^2  ~~~~ \text{subject to} ~~~~Q(\bold{x})=\bold{0}.
\end{eqnarray} 
\end{subequations} 
The augmented Lagrangian function for the problem is \citep{Nocedal_2006_NOO}
\begin{equation} \label{app2}
\mathcal{L}_A(\bold{x},\bold{v}) =\|P(\bold{x})\|_2^2 +  \bold{v}^T Q(\bold{x}) + \frac{\lambda}{2}\|Q(\bold{x})\|_2^2.
\end{equation}
The problem \ref{app2} can be written in a more compact form by introducing the scaled dual variable $\bar{\bold{q}}=-\frac{\bold{v}}{\lambda}$ and adding and subtracting the term $\frac{\lambda}{2}\|\bar{\bold{q}}\|_2^2$ to the augmented Lagrangian \ref{app2}. In this case, we arrive at the following scaled-form of the method of multipliers:
\begin{subequations}
\label{app3}
\begin{eqnarray}
\mathcal{L}_A(\bold{x},\bar{\bold{q}}) &=&\|P(\bold{x})\|_2^2 - \lambda \bold{\bar{\bold{q}}}^T Q(\bold{x}) + \frac{\lambda}{2}\|Q(\bold{x})\|_2^2+\frac{\lambda}{2}\|\bar{\bold{q}}\|_2^2-\frac{\lambda}{2}\|\bar{\bold{q}}\|_2^2\\
&=&\|P(\bold{x})\|_2^2 + \frac{\lambda}{2} \|Q(\bold{x})-\bar{\bold{q}}\|_2^2 - \frac{\lambda}{2}\|\bar{\bold{q}}\|_2^2.
\end{eqnarray} 
\end{subequations}

\bibliographystyle{seg}

\newcommand{\SortNoop}[1]{}

\end{document}